\documentclass{article}

\usepackage{amsmath}
\usepackage{amssymb}
\usepackage{amsfonts}

\begin{document}

\renewcommand{\theequation}{\arabic{section}.\arabic{equation}}

\begin{center}
{\Large Rescaling approach for a stochastic population dynamics equation
perturbed by a linear multiplicative Gaussian noise}

\bigskip

{\large Gabriela Marinoschi}

\bigskip

\textquotedblleft Gheorghe Mihoc-Caius Iacob\textquotedblright\ Institute of
Mathematical Statistics and Applied Mathematics of the Romanian Academy,
Calea 13 Septembrie 13, 050711 Bucharest, Romania

E-mail address: gabriela.marinoschi@acad.ro
\end{center}

\bigskip

Abstract. We are concerned with a nonlinear nonautonomous model represented
by an equation describing the dynamics of an age-structured population
diffusing in a space habitat $O,$ governed by local Lipschitz vital factors
and by a stochastic behavior of the demographic rates possibly representing
emigration, immigration and fortuitous mortality. The model is completed by
a random initial condition, a flux type boundary conditions on $\partial O$
with a random jump in the population density and a nonlocal nonlinear
boundary condition given at age zero. The stochastic influence is expressed
by a linear multiplicative Gaussian noise perturbation in the equation. The
main result proves that the stochastic model is well-posed, the solution
being in the class of path-wise continuous functions and satisfying some
particular regularities with respect to the age and space. The approach is
based on a rescaling transformation of the stochastic equation into a random
deterministic time dependent hyperbolic-parabolic equation with local
Lipschitz nonlinearities. The existence and uniqueness of a strong solution
to the random deterministic equation is proved by combined semigroup,
variational and approximation techniques. The information given by these
results is transported back via the rescaling transformation towards the
stochastic equation and enables the proof of its well-posedness.

\bigskip

Keywords: stochastic nonlinear equation, noise induced phenomena,
multiplicative noise, Brownian motion, population dynamics, semigroup
approach

2010 Mathematics Subject Classification: 35R60, 60H15, 92D25, 35Q92

\bigskip

\section{Problem presentation}

\setcounter{equation}{0}

This paper addresses the well-posedness of a nonlinear stochastic population
dynamics equation describing the time and age evolution of a population of
density $p,$ in an space habitat $O,$ governed by nonlinear vital factors,
as natality and mortality and environmental particularities and influenced
by a linear multiplicative noise perturbation. The equation reads 
\begin{eqnarray}
&&dp(t,a,x)+p_{a}(t,a,x)dt-\Delta p(t,a,x)dt+\mu _{S}(t,a,x;U(p))p(t,a,x)dt
\label{1} \\
&=&p(t,a,x)dW(t,a,x),\mbox{ in }(0,T)\times (0,a^{+})\times O.  \notag
\end{eqnarray}%
It is completed by two boundary conditions, the first of Robin type on the
boundary of $O,$ and the second at $a=0,$ and by the initial condition (\ref%
{4}), below: 
\begin{equation}
-\nabla p(t,a,x)\cdot \nu =\alpha _{0}(t,a,x)p(t,a,x)+k_{0}(t,a,x),\mbox{ on 
}(0,T)\times (0,a^{+})\times \partial O,  \label{2}
\end{equation}%
\begin{equation}
p(t,0,x)=\int_{0}^{a^{+}}m_{0}(a,x;U(p))p(t,a,x)da,\mbox{ in }(0,T)\times O,
\label{3}
\end{equation}%
\begin{equation}
p(0,a,x)=p_{0}(a,x),\mbox{ in }(0,a^{+})\times O.  \label{4}
\end{equation}%
In these equations, $t$ is the time running in $(0,T),$ $a$ is the age
belonging to $(0,a^{+}),$ $a^{+}$ is the maximum age life and $x$ is the
space variable in $O$ which is an open bounded domain of $\mathbb{R}^{d}$ ($%
d=1,2,3).$ The Laplacian $\Delta $ and the gradient $\nabla $ refer\ only to
the space variable and $\nu $ is the unit outward normal vector to the
boundary of $O.$ Moreover, $\mu _{S},$ called supplementary or additional
mortality, is the mortality rate due to other causes than reaching the
maximum life $a^{+}$, and $m_{0}$ is the fertility rate. The vital rates are
allowed to depend nonlinearly on $p,$ by the variable 
\begin{equation}
U(p)=\int_{0}^{a^{+}}\int_{O_{U}}\gamma (a,x)p(t,a,x)dxda,  \label{U}
\end{equation}%
where $O_{U}$ is a subset of $O.$ This dependence means that the total
population lying in the environment $O_{U}$ may produce a perturbation of
the vital rates according to the weight factor $\gamma $ varying with
respect to age and space. It is relevant \ to assume that the dependence of $%
\mu _{S}$ and $m_{0}$ on the variable $U$ is locally Lipschitz (see e.g., 
\cite{cuiama-1} and \cite{cuiama-2}).

The boundary condition (\ref{3}) written for $a=0$ is the well-known birth
equation in population dynamics. The boundary condition (\ref{2}) expresses
a change of population living in the habitat $O$ with the outer environment,
supplemented by a possible jump in the population density on the boundary.
Other types of conditions, indicating a hostile boundary or a closed
habitat, can be considered by assuming homogeneous Dirichlet or Neumann
boundary conditions, respectively.

\noindent We note that the population dynamics equation with age-structure (%
\ref{1}) normally includes also, on the left-hand side, a term $\mu
_{0}(a)p, $ where $\mu _{0}(a)$ is the natural mortality due to reaching the
maximum life age. Since the natural assumptions for $\mu _{0}$ indicates
that this is a $L_{loc}^{1}$ function in $(0,a^{+})$, a standard treatment
is to replace\ $p$ by $p\exp \left( -\int_{0}^{a}\mu _{0}(s)ds\right) .$ In
this way the term $\mu _{0}(a)p$ is cancelled from the equation and so,
without loss of generality, the equation reduces to the so-called normalized
equation (\ref{1}).

Now, let us pass to the stochastic context. A deterministic model (with zero
on the right-hand side of (\ref{1})) cannot reproduce or explain the effects
of random fluctuations which come from the intrinsic stochastic nature of
open systems. Random effects may be also induced by the interplay between
the behavior of natural systems and random fluctuations generated by the
environment. The presence of noise produced by this interaction determines
an increase of the complexity of the system evolution which can
substantially drift apart from its known deterministic feature. Moreover,
demographic events basically represented by statistical averages lead to a
weak determinism and so, in order to describe all these, a pure stochastic
contribution should be taken into account in the equation.

Let ($\Omega ,\mathcal{F},\mathbb{P})$ be a probability space, with the
natural filtration $\{\mathcal{F}_{t}\}_{t\geq 0}$ and let $W$ be a
stochastic Gaussian process of the form 
\begin{equation}
W(t,a,x)=\sum_{j=1}^{N}\mu _{j}(a,x)\beta _{j}(t),  \label{5}
\end{equation}%
where $\{\beta _{j}\}_{j=1}^{N}$ is an independent system of real-valued
Brownian motions and $\mu _{j}$ are regular functions. Thus, in relation
with equation (\ref{1}), $W$ mimics a random fluctuation in the interaction
between the population and the environment which can be due to a possible
immigration, emigration or incidental stochastic mortality.

Let us note that $W$ can be taken as well as an infinite series, like e.g.,
in \cite{vb-mr-archive} under certain convergence conditions for the series
of the square coefficients. As usually, the argument $\omega \in \Omega $
will be not explicitly specified in the list of the function arguments.

The deterministic nonlinear model (for $W=0)$ described by equations of type
(\ref{1})-(\ref{4}) made the subject of a large amount of papers in the
literature. A synthetic presentation of the most important achievements can
be found in the monograph \cite{Iannelli-Milner} and in the references
therein, and in relation with the deterministic nonlinear model with locally
Lipschitz nonlinearities for the vital rates, in \cite{cuiama-1} and \cite%
{cuiama-2}.

The autonomous stochastic linear model of type (\ref{1})-(\ref{4}),
characterized by $\gamma =0$ and time independent coefficients $\mu _{S},$ $%
\alpha _{0},$ $k_{0},$ can be analyzed by rewriting this system in an
operatorial form and treating it by a semigroup formulation in the $L^{2}$%
-approach, as e.g., in \cite{DaPrato-Zabcyk}, or \cite{prevot-MR}. A
path-wise continuous solution to the linear autonomous stochastic equation
can be proved, if $\mu _{S}$ is globally Lipschitz continuous. Under a
supplementary condition on the operator, the well-posedness may follow in
the case of a local Lipschitz term $\mu _{S}$ for the stochastic equation
with an additive noise (see e.g., \cite{DaPrato-Zabcyk}, chapter 7). We also
refer to the paper \cite{vb-mr-2015} in which the existence for stochastic
equations with a linear multiplicative noise, with a general nonlinear
monotone, demicontinuous and coercive time dependent operator between two
dual spaces is provided. However, these results are not directly applicable
in our case because the problem is not parabolic-like as in \cite{vb-mr-2015}
and the nonlinearities are not globally Lipschitz. As far as we know, the
stochastic equation (\ref{1}) with $m_{0}$ and $\mu _{S}$ local Lipschitz
has been not addressed in the literature.

The proof we provide begins by applying to our problem a rescaling
transformation. More exactly, by a suitable function transformation for $p$,
system (\ref{1})-(\ref{4}) is transformed into the random deterministic one,
in the unknown $y,$ see (\ref{26})-(\ref{29}) in Section 2. This is a
nonlinear time dependent hyperbolic-parabolic system with local Lipschitz
nonlinearities and it cannot be fitted in any functional framework for which
general existence theorems can be applied. The proof of the solution
existence requires a long and technical approach which is split in many
intermediate results, beginning with the study of the well-posedness of a
generic hyperbolic-parabolic problem with globally Lipschitz nonlinearities,
in Section 3, Proposition 3.2. This proof is led by means of combined
semigroup, variational\ and approximation techniques. Relying on this result
and by using two types of regularizations, one for the time coefficients and
the other for the operator acting in the equation, the existence and
uniqueness of the solution to the random system is given in Theorem 4.1, for
all $\omega \in \Omega $. Much effort is done to get estimates for the
solution to the intermediate problems in order to ensure the strong
convergence in the approximated equations. In addition, some space and age
regularity and the existence of a strong solution for the rescaled equation
are obtained in Corollary 4.2. All information provided by the solution to
the random system is fructified, while going back via the rescaling
transformation, into the proof of well-posedness of the stochastic system,
in Theorem 4.3.

\medskip

\paragraph{Notation.}

For $1\leq p\leq \infty $, $L^{p}(O)$ is the space of all $p$-integrable
real valued functions with the norm $\left\Vert \cdot \right\Vert
_{L^{p}(O)} $ and $L^{q}(0,T;L^{p}(O))$ denotes the space of measurable
functions $u:[0,T]\rightarrow L^{p}(O)$ such that $t\rightarrow \left\Vert
u(t)\right\Vert _{L^{p}(O)}$ belongs to $L^{q}(0,T).$ $C([0,T],L^{p}(O))$ is
the continuous $L^{p}$-valued functions with the supremum norm in $t.$ As
usually, $W^{1,p}(O)$ is the classical Sobolev space, i.e., $%
W^{1,p}(O)=\{u\in L^{p}(O);$ $\nabla u\in L^{p}(O)\}$ and $%
H^{1}(O)=W^{1,2}(O).$ The scalar product and the norm in a Hilbert space $X$
are denoted by $\left( \cdot ,\cdot \right) _{X}$ and $\left\Vert \cdot
\right\Vert _{X},$ respectively. In particular, $\left\Vert \cdot
\right\Vert _{\infty }$ indicates the norm of functions belonging to $%
L^{\infty }(0,T)\times (0,a^{+})\times O)$ or $L^{\infty }(0,T)\times
(0,a^{+})\times \partial O).$

If no confusion can be done, some function arguments will be not specified
in the integrands. $C,$ $C_{i},$ $c_{i},$ $i=0,1,2,...$ will stand for
several constants that may change in the computations from line to line.
Moreover, we shall denote 
\begin{equation*}
H=L^{2}(O),\mbox{ \ }V=H^{1}(O),\mbox{ \ }V^{\prime }=(H^{1}(O))^{\prime },
\end{equation*}%
\begin{equation*}
\mathcal{H}=L^{2}(0,a^{+};H),\mbox{ \ }\mathcal{V}=L^{2}(0,a^{+};V),\mbox{ \ 
}\mathcal{V}^{\prime }=L^{2}(0,a^{+};V^{\prime }).
\end{equation*}%
where $V^{\prime }$ is the dual of $V,$ and $\mathcal{V}^{\prime }$ is the
dual of $\mathcal{V}$. We note that $V\subset H\subset V^{\prime }$ and $%
\mathcal{V\subset H\subset V}^{\prime }$ with compact injections.

\section{Preliminaries}

\setcounter{equation}{0}

We specify the hypotheses which will be in effect in this work (see e.g., 
\cite{cuiama-1}).

We assume, as relevant in population dynamics, that $\mu _{S}(t,a,x;r)$ and $%
m_{0}(a,x;r)$ are local Lipschitz functions on $\mathbb{R}$ in the variable $%
r,$ uniformly with respect to $t,a,x,$ i.e., for any $R>0,$ there exists $%
L_{\mu _{S}}(R)$ and $L_{m_{0}}(R)$ such that 
\begin{eqnarray}
\left\vert \mu _{S}(t,a,x;r)-\mu _{S}(t,a,x;\overline{r})\right\vert &\leq
&L_{\mu _{S}}(R)\left\vert r-\overline{r}\right\vert ,  \label{9-Lip} \\
\left\vert m_{0}(a,x;r)-m_{0}(a,x;\overline{r})\right\vert &\leq
&L_{m_{0}}(R)\left\vert r-\overline{r}\right\vert ,  \notag
\end{eqnarray}%
whenever $\left\vert r\right\vert \leq R$ and $\left\vert \overline{r}%
\right\vert \leq R$. Moreover, 
\begin{eqnarray}
\mu _{S}(\cdot ,\cdot ,\cdot ;r) &\in &L^{\infty }((0,T)\times
(0,a^{+})\times O),\mbox{ for all }r\in \mathbb{R},  \label{8} \\
0 &\leq &\mu _{S}(t,a,x;r)\leq \mu _{\infty }\mbox{ a.e. in }(0,T)\times
(0,a^{+})\times O,\mbox{ for all }r\in \mathbb{R},  \notag
\end{eqnarray}%
\begin{eqnarray}
m_{0}(\cdot ,\cdot ;r) &\in &L^{\infty }((0,a^{+})\times O),\mbox{ for all }%
r\in \mathbb{R},  \label{9} \\
0 &\leq &m_{0}\leq m_{0\infty }\mbox{ a.e. in }(0,a^{+})\times O,  \notag
\end{eqnarray}%
\begin{equation*}
\gamma \in L^{\infty }((0,a^{+})\times O),\mbox{ \ }0\leq \gamma (a,x)\leq
\gamma _{\infty }\mbox{ a.e. in }(0,a^{+})\times O.
\end{equation*}%
We still assume that%
\begin{equation}
\alpha _{0}\in L^{\infty }((0,T)\times (0,a^{+})\times \partial O),\mbox{ \ }%
\alpha _{0}\geq 0\mbox{ a.e. in }(0,T)\times (0,a^{+})\times O.  \label{10}
\end{equation}%
Moreover, $k_{0}$ and $p_{0}$ are random functions$,$ such that%
\begin{equation}
k_{0}\in L^{2}((0,T)\times (0,a^{+})\times \partial O),\mbox{ }\mathbb{P}%
\mbox{-a.s.,}  \label{10-0}
\end{equation}%
\begin{equation}
p_{0}\in L^{2}(0,a^{+};L^{2}(O)),\mbox{ }\mathbb{P}\mbox{-a.s.,}  \label{11}
\end{equation}%
\begin{equation}
p_{0}(\cdot ,a,x)\mbox{ is measurable with respect to }\mathcal{F}_{0},\mbox{
a.a. }(a,x),  \label{11-0}
\end{equation}%
\begin{equation*}
k_{0}(\cdot ,t,a,x)\mbox{ is }\mathcal{F}_{t}\mbox{-adapted},\mbox{ a.a. }%
(t,a,x).
\end{equation*}%
Finally, $\beta _{j}\in C([0,T];\mathbb{R}),$ \ $\beta _{j}(0)=0,$ $%
j=1,...,N,$ and assume that%
\begin{equation}
\mu _{j}\in C^{2}([0,a^{+}]\times \overline{O}),\mbox{ \ }\nabla \mu
_{j}\cdot \nu =0\mbox{ on }(0,a^{+})\times \partial O,\mbox{ \ }j=1,...,N.
\label{6}
\end{equation}%
In particular, for a.a. $\omega \in \Omega ,$ the mapping $%
(t,a,x)\rightarrow W(t,a,x)(\omega )$ is continuous and the process $%
\{W(t,a,x)\}_{t\geq 0}$ is real-valued $\mathcal{F}_{t}$-adapted. As
usually, we shall not specify the variable $\omega $ in all random functions
that occur.

\medskip

\noindent \textbf{Definition 2.1. }A process $p:[0,T]\times \Omega
\rightarrow \mathcal{H}$ is called a solution to (\ref{1})-(\ref{4}) if it
is an $\mathcal{F}_{t}$-adapted process, $t\geq 0,$ 
\begin{equation}
p\in C([0,T];\mathcal{H})\cap L^{2}(0,T;\mathcal{V})\cap
C([0,a^{+}];L^{2}(0,T;H)),\mbox{ }\mathbb{P}\mbox{-a.s},\mbox{ }  \label{20}
\end{equation}%
and 
\begin{eqnarray}
&&\left( p(t),\psi \right) _{\mathcal{H}}+\int_{0}^{t}\int_{O}p(\tau
,a^{+},x)\psi (a^{+},x)dxd\tau -\int_{0}^{t}\int_{0}^{a^{+}}\int_{O}p\psi
_{a}dxdad\tau  \label{21} \\
&&-\int_{0}^{t}\int_{0}^{a^{+}}\int_{O}m_{0}(a,x;U(p))p\psi (0,x)dxdad\tau 
\notag \\
&&+\int_{0}^{t}\int_{0}^{a^{+}}\int_{O}(\nabla p\cdot \nabla \psi +\mu
_{S}(\tau ,a,x;U(p))p\psi )dxdad\tau
+\int_{0}^{t}\int_{0}^{a^{+}}\int_{\partial O}(\alpha _{0}p+k_{0})\psi
d\sigma dad\tau  \notag \\
&=&\left( p_{0},\psi \right) _{\mathcal{H}}+\int_{0}^{t}\left( p(\tau
)dW(\tau ),\psi \right) _{\mathcal{H}},\mbox{ }\mathbb{P}\mbox{-a.s., for
all }\psi \in \mathcal{V}\mbox{, with }\psi _{a}\in \mathcal{V}^{\prime }. 
\notag
\end{eqnarray}

\medskip

\noindent We specify that since $p\in C([0,T];\mathcal{H})$ $\mathbb{P}$-a.s$%
,$ the It\^{o} integral 
\begin{equation}
\int_{0}^{t}\left( p(\tau )dW(\tau ),\psi \right) _{\mathcal{H}%
}=\sum_{j=1}^{N}\int_{0}^{t}\left( \int_{0}^{a^{+}}\int_{O}\mu _{j}(a,x)\psi
(a,x)p(\tau ,a,x)dxda\right) d\beta _{j}(\tau )  \label{22}
\end{equation}%
is well defined.

We begin by transforming equation (\ref{1}), using for this the rescaling
formula%
\begin{equation}
p(t,a,x)=e^{W(t,a,x)}y(t,a,x),\mbox{ for }t\geq 0.  \label{23}
\end{equation}%
In the calculations implied by this transformation we use the It\^{o}'s
relation%
\begin{equation}
de^{W}=e^{W}dW+\mu e^{W}dt  \label{23-0}
\end{equation}%
where, 
\begin{equation}
\mu (a,x)=\frac{1}{2}\sum_{j=1}^{N}\mu _{j}^{2}(a,x).  \label{25}
\end{equation}%
Then, (\ref{23})-(\ref{25}) imply that%
\begin{equation}
dp=e^{W}dy+e^{W}ydW+\mu e^{W}ydt.  \label{24}
\end{equation}%
After plugging (\ref{23}) in (\ref{1})-(\ref{4}) and performing some
calculations by expanding the terms $\Delta (e^{W}y)$ and $\nabla (e^{W}y)$
we deduce the system 
\begin{equation}
y_{t}+y_{a}-\Delta y+g_{1}(t,a,x)y+g_{2}(t,a,x)\cdot \nabla y+\mu
_{S}(t,a,x;U(e^{W}y))y=0,\mbox{ in }(0,T)\times (0,a^{+})\times O,
\label{26}
\end{equation}%
\begin{equation}
-\nabla y\cdot \nu =\alpha (t,a,x)y+k(t,a,x),\mbox{ in }(0,T)\times
(0,a^{+})\times \partial O,  \label{27}
\end{equation}%
\begin{equation}
y(t,0,x)=\int_{0}^{a^{+}}m(t,a,x;U(e^{W}y))y(t,a,x)da,\mbox{ in }(0,T)\times
O,  \label{28}
\end{equation}%
\begin{equation}
y(0,a,x)=y_{0}(a,x)=p_{0}(a,x),\mbox{ in }(0,a^{+})\times O,  \label{29}
\end{equation}%
where 
\begin{eqnarray}
g_{1} &=&W_{a}-\Delta W-\left\vert \nabla W\right\vert ^{2}+\mu ,\mbox{ \ }%
g_{2}=-2\nabla W,  \label{30} \\[0.01in]
\alpha &=&\alpha _{0}+\nabla W\cdot \nu =\alpha _{0},\mbox{ \ }k=k_{0}e^{-W},
\notag \\
m(t,a,x;r) &=&m_{0}(a,x;r)e^{W(t,a,x)-W(t,0,x)}.  \notag
\end{eqnarray}%
The functions $g_{1},$ $g_{2},$ $\alpha $ and $k$ depend on $t,$ $a,$ $x,$
and $\alpha \geq 0$ by (\ref{10}) and (\ref{6}) and obviously, the functions 
$\mu _{S}$ and $m$ are locally Lipschitz continuous with respect to the
fourth variable, with the Lipschitz constants $L_{\mu _{S}}(R)$ and $%
L_{m}(R)=L_{m_{0}}(R)c_{W_{0}},$ where 
\begin{equation}
c_{W_{0}}=\left\Vert e^{W(\cdot ,\cdot ,\cdot )-W(\cdot ,0,\cdot
)}\right\Vert _{\infty }.  \label{cW}
\end{equation}

On behalf of the hypotheses (\ref{9-Lip})-(\ref{6}) we deduce%
\begin{eqnarray}
g_{1} &\in &C([0,T];C^{1}([0,a^{+}]\times C(\overline{O})),\mbox{ }g_{2}\in
C([0,T];C^{2}[0,a^{+}]\times C^{1}(\overline{O})),  \label{31} \\
\alpha &=&\alpha _{0}\in L^{\infty }((0,T)\times (0,a^{+})\times \partial
O)),\mbox{ \ }k\in L^{2}((0,T)\times (0,a^{+})\times \partial O).  \notag
\end{eqnarray}%
It is obvious that (\ref{26})-(\ref{29}) is deterministic but random.

\medskip

\noindent \textbf{Definition 2.2. }A solution $y:[0,T]\times \Omega
\rightarrow \mathcal{H}$  to (\ref{26})-(\ref{29}) is an $\mathcal{F}_{t}$%
\textit{-}adapted process such that 
\begin{equation}
y\in C([0,T];\mathcal{H})\cap L^{2}(0,T;\mathcal{V})\cap
C([0,a^{+}];L^{2}(0,T;H)),\mbox{ }\mathbb{P}\mbox{-a.s.},  \label{32}
\end{equation}%
and 
\begin{eqnarray}
&&-\int_{0}^{T}\int_{0}^{a^{+}}\int_{O}y\psi
_{t}dxdadt-\int_{0}^{a^{+}}\int_{O}y_{0}\psi (0,a,x)dxda  \label{33} \\
&&+\int_{0}^{T}\int_{O}y(t,a^{+},x)\psi
(t,a^{+},x)dxdt-\int_{0}^{T}\int_{0}^{a^{+}}\int_{O}y\psi _{a}dxdadt  \notag
\\
&&-\int_{0}^{T}\int_{O}\left( \int_{0}^{a^{+}}m(t,a,x;U(e^{W}y))da\right)
\psi (t,0,x)dxdt+\int_{0}^{T}\int_{0}^{a^{+}}\int_{\partial O}(\alpha y\psi
+k\psi )d\sigma dadt  \notag \\
&&+\int_{0}^{T}\int_{0}^{a^{+}}\int_{O}\left( \nabla y\cdot \nabla \psi
+yg_{1}\psi +\psi g_{2}\cdot \nabla y+\mu _{S}(t,a,x;U(e^{W}y))y\psi \right)
dxdadt=0,\mbox{ }\mathbb{P}\mbox{-a.s.,}  \notag
\end{eqnarray}%
for all $\psi \in W^{1,2}(0,T;\mathcal{H})\cap L^{2}(0,T;\mathcal{V)}$, with 
$\psi _{a}\in L^{2}(0,T;\mathcal{V}^{\prime })$ and $\psi (T,a,x)=0.$

\section{Intermediate results}

\setcounter{equation}{0}

As we shall see, due to the local Lipschitz properties of $m$ and $\mu _{S}$
the proof of the existence of the solution to the deterministic random
system is very long and technical. For making the arguments more
understandable we shall split it in many parts. We begin with an
intermediate result for a generic deterministic hyperbolic-parabolic time
dependent system with global Lipschitz nonlinearities.

We introduce the problem 
\begin{equation}
Y_{t}+Y_{a}-\Delta Y+f_{1}(t,a,x)Y+f_{2}(t,a,x)\cdot \nabla
Y+E_{1}(t,a,x;Y)=f,\mbox{ in }(0,T)\times (0,a^{+})\times O,  \label{44}
\end{equation}%
\begin{equation}
-\nabla Y\cdot \nu =Yf_{\Gamma }+f_{\Gamma }^{0},\mbox{ in }(0,T)\times
(0,a^{+})\times \partial O,  \label{45}
\end{equation}%
\begin{equation}
Y(t,0,x)=\int_{0}^{a^{+}}E_{2}(t,a,x;Y)da,\mbox{ in }(0,T)\times O,
\label{46}
\end{equation}%
\begin{equation}
Y(0,a,x)=Y_{0}(a,x),\mbox{ in }(0,a^{+})\times O,  \label{47}
\end{equation}%
where 
\begin{eqnarray}
f_{1} &\in &C^{1}([0,T];C^{1}([0,a^{+}]\times C(\overline{O})),\mbox{ }
\label{48} \\
f_{2} &\in &C^{1}([0,T];C^{2}[0,a^{+}]\times C^{1}(\overline{O})),\mbox{ }%
f_{2}\cdot \nu =0\mbox{ on }(0,T)\times (0,a^{+})\times \partial O,  \notag
\\
f_{\Gamma } &\in &C^{1}([0,T];L^{\infty }((0,a^{+})\times \partial O)),\mbox{
\ }f_{\Gamma }(t,a,x)\geq 0,\mbox{ a.e. in }(0,T)\times (0,a^{+})\times
\partial O,  \notag \\
f_{\Gamma }^{0} &\in &L^{2}(0,T;L^{2}((0,a^{+})\times \partial O)),\mbox{ }%
f\in L^{2}(0,T;\mathcal{H}),\mbox{ }Y_{0}\in \mathcal{H}.  \notag
\end{eqnarray}%
Here, $E_{i}:(0,T)\times (0,a^{+})\times O\times \mathcal{H\rightarrow H}$, $%
i=1,2,$ and both operators are globally Lipschitz on $\mathcal{H}$,
uniformly for $(t,a,x)\in (0,T)\times (0,a^{+})\times O,$ i.e., there exist $%
L_{i}>0,$ $i=1,2,$ such that for any $v,$ $\overline{v}\in \mathcal{H}$ we
have 
\begin{equation}
\left\Vert E_{1i}(t,\cdot ,\cdot ;v)-E_{i}(t,\cdot ,\cdot ;\overline{v}%
)\right\Vert _{\mathcal{H}}\leq L_{i}\left\Vert v-\overline{v}\right\Vert _{%
\mathcal{H}},  \label{Li}
\end{equation}%
for any $t\in (0,T).$ Moreover, 
\begin{equation}
\left\Vert E_{1}(t,\cdot ,\cdot ;v)\right\Vert _{\mathcal{H}}\leq \mu
_{\infty }\left\Vert v\right\Vert _{\mathcal{H}},\mbox{ \ }\left\Vert
E_{2}(t,\cdot ,\cdot ;v)\right\Vert _{\mathcal{H}}\leq m_{\infty }\left\Vert
v\right\Vert _{\mathcal{H}},\mbox{ for all }v\in \mathcal{H},\mbox{ }
\label{cE}
\end{equation}%
uniformly with respect to $(t,a,x)$ where $\mu _{\infty }$ and $m_{\infty }$
are precisely given by 
\begin{eqnarray*}
\mu _{\infty } &:&=\sup_{(t,a,x,r)\in (0,T)\times (0,a^{+})\times O\times 
\mathbb{R}}\left\vert \mu _{S}(t,a,x;r)\right\vert ,\mbox{ \ } \\
m_{\infty } &:&=\sup_{(t,a,x,r)\in (0,T)\times (0,a^{+})\times O\times 
\mathbb{R}}\left\vert m(t,a,x;r)\right\vert =c_{W_{0}}m_{0\infty }.
\end{eqnarray*}

\medskip

\noindent \textbf{Definition 3.1. }A solution $Y$ to (\ref{44})-(\ref{47})
is a function 
\begin{equation}
Y\in C([0,T];\mathcal{H})\cap C([0,a^{+}];L^{2}(0,T;H))\cap L^{2}(0,T;%
\mathcal{V})  \label{48-00}
\end{equation}%
which satisfies the equation 
\begin{eqnarray}
&&-\int_{0}^{T}\int_{0}^{a^{+}}\int_{O}Y\psi
_{t}dxdadt-\int_{0}^{a^{+}}\int_{O}Y_{0}\psi
(0,a,x)dxda-\int_{0}^{T}\int_{0}^{a^{+}}\int_{O}Y\psi _{a}dxdadt
\label{48-0} \\
&&+\int_{0}^{T}\int_{O}Y(t,a^{+},x)\psi
(t,a^{+},x)dxdt-\int_{0}^{T}\int_{O}\left(
\int_{0}^{a^{+}}E_{2}(t,a,x;Y)da\right) \psi (t,0,x)dxdt  \notag \\
&&+\int_{0}^{T}\int_{0}^{a^{+}}\int_{O}(\nabla Y\cdot \nabla \psi
+f_{1}Y\psi +\psi f_{2}\cdot \nabla Y)dxdadt+  \notag \\
&&\int_{0}^{T}\int_{0}^{a^{+}}\int_{\partial O}(f_{\Gamma }Y\psi +f_{\Gamma
}^{0}\psi )d\sigma
dadt+\int_{0}^{T}\int_{0}^{a^{+}}\int_{O}E_{1}(t,a,x;Y)\psi
dxdadt=\int_{0}^{T}\int_{0}^{a^{+}}\int_{O}f\psi dxdadt,\mbox{ }  \notag
\end{eqnarray}%
for all $\psi \in W^{1,2}(0,T;\mathcal{H})\cap L^{2}(0,T;\mathcal{V})$, with 
$\psi _{a}\in L^{2}(0,T;\mathcal{V}^{\prime })$\ and $\psi (T,a,x)=0.$\ 

\medskip

\noindent \textbf{Proposition 3.2. }\textit{Under the assumptions} (\ref{48}%
), \textit{problem}\textbf{\ }(\ref{44})-(\ref{47}) \textit{has a unique
solution, which satisfies the estimate}%
\begin{eqnarray}
&&\left\Vert Y(t)\right\Vert _{\mathcal{H}}^{2}+\int_{0}^{t}%
\int_{0}^{a}Y^{2}(\tau ,a,x)dxd\tau +\int_{0}^{t}\left\Vert Y(\tau
)\right\Vert _{\mathcal{V}}^{2}d\tau  \label{48-1} \\
&\leq &\overline{\overline{C}}:=C_{est}\left( \left\Vert Y_{0}\right\Vert _{%
\mathcal{H}}^{2}+\int_{0}^{t}\left\Vert f_{\Gamma }^{0}(\tau )\right\Vert
_{L^{2}(0,a^{+};L^{2}(\partial O))}^{2}d\tau +\int_{0}^{t}\left\Vert f(\tau
)\right\Vert _{\mathcal{H}}^{2}d\tau \right) ,  \notag
\end{eqnarray}%
\textit{for all }$t\in \lbrack 0,T]$ and $a\in \lbrack 0,a^{+}],$\textit{\
where}%
\begin{equation}
C_{est}=c_{0}e^{c_{1}(1+\left\Vert f_{1}\right\Vert _{\infty }+\left\Vert
f_{2}\right\Vert _{\infty }^{2}+a^{+}m_{0\infty }^{2}c_{W_{0}}^{2}+\mu
_{\infty }^{2})T},  \label{Cest}
\end{equation}%
$c_{0}$\textit{\ and }$c_{1}$\textit{\ are positive numbers, }$c_{0}$ 
\textit{depending on the domain and space dimension.} \textit{Moreover, for
two solutions }$Y$ \textit{and }$\overline{Y}$\textit{\ corresponding to }$%
\{Y_{0},$\textit{\ }$f_{1}^{1},$ $f_{2}^{1},$ $f_{\Gamma }^{1},$\textit{\ }$%
f_{\Gamma }^{01},$\textit{\ }$f^{1}\}$\textit{\ and }$\{Y_{0}^{2},f_{1}^{2},$
$f_{2}^{2},$ \textit{\ }$f_{\Gamma }^{2},$\textit{\ }$f^{02},$ $f^{2}\},$%
\textit{\ respectively, we have}%
\begin{eqnarray}
&&\left\Vert (Y_{1}-Y_{2})(t)\right\Vert _{\mathcal{H}}^{2}+\int_{0}^{t}%
\int_{O}(Y_{1}-Y_{2})^{2}(\tau ,a,x)dxd\tau +\int_{0}^{t}\left\Vert
(Y_{1}-Y_{2})(\tau )\right\Vert _{\mathcal{V}}^{2}d\tau  \label{700} \\
&\leq &c_{0}e^{c_{1}(1+\left\Vert f_{1}^{2}\right\Vert _{\infty }+\left\Vert
\nabla \cdot f_{2}^{2}\right\Vert _{\infty
}+L_{1}^{2}+a^{+}L_{2}^{2})T}\times \left( \left\Vert
Y_{0}^{1}-Y_{0}^{2}\right\Vert _{\mathcal{H}}^{2}+\int_{0}^{t}\left\Vert
(f^{1}-f^{2})(\tau )\right\Vert _{\mathcal{H}}^{2}d\tau \right.  \notag \\
&&\left. +\overline{\overline{C}}\left( \left\Vert
f_{1}^{1}-f_{1}^{2}\right\Vert _{\infty }^{2}+\left\Vert
f_{2}^{1}-f_{2}^{2}\right\Vert _{\infty }^{2}+\left\Vert f_{\Gamma
}^{1}-f_{\Gamma }^{2}\right\Vert _{\infty }^{2}\right)
+\int_{0}^{t}\left\Vert (f_{\Gamma }^{01}-f_{\Gamma }^{02})(\tau
)\right\Vert _{L^{2}(0,a^{+};L^{2}(\partial O))}^{2}d\tau \right) ,  \notag
\end{eqnarray}%
\textit{for all }$t\in \lbrack 0,T]$ \textit{and} $a\in \lbrack 0,a^{+}].$%
\textit{\ }

\medskip

\noindent \textbf{Proof. }The proof is done in 4 steps.

\noindent \textbf{Step 1.} Let us consider $E_{1}=E_{2}=f_{\Gamma }^{0}=0.$
For all $t\in \lbrack 0,T]$ we introduce the time dependent operator $%
A_{0}(t):\mathcal{V\rightarrow V}^{\prime }$ by 
\begin{eqnarray*}
\left\langle A_{0}(t)v,\psi \right\rangle _{\mathcal{V}^{\prime },\mathcal{V}%
} &=&\left\langle v_{a},\psi \right\rangle _{\mathcal{V}^{\prime },\mathcal{V%
}}+\int_{0}^{a^{+}}\int_{O}\nabla v\cdot \nabla \psi
dxda+\int_{0}^{a^{+}}\int_{\partial O}vf_{\Gamma }(t,a,x)\psi d\sigma da \\
&&+\int_{0}^{a^{+}}\int_{O}\left( f_{1}(t,a,x)v+\nabla v\cdot
f_{2}(t,a,x)\right) \psi dxda,\mbox{ for all }\psi \in \mathcal{V}.
\end{eqnarray*}%
We specify that $\left\langle \cdot ,\cdot \right\rangle _{\mathcal{V}%
^{\prime },\mathcal{V}}$ is the pairing between the dual spaces $\mathcal{V}%
^{\prime }$ and $\mathcal{V},$ defined as 
\begin{equation*}
\left\langle \phi ,\psi \right\rangle _{\mathcal{V}^{\prime },\mathcal{V}%
}=\int_{0}^{a^{+}}\left\langle \phi (a),\psi (a)\right\rangle _{V^{\prime
},V}da,\mbox{ for }\phi \in \mathcal{V}^{\prime },\mbox{ }\psi \in \mathcal{V%
}.
\end{equation*}%
Next, we define the restriction $A(t):D(A(t))\subset \mathcal{H}\rightarrow 
\mathcal{H}$, where 
\begin{equation*}
D(A(t))=\{v\in \mathcal{V};v_{a}\in \mathcal{V}^{\prime },v(0,x)=0,A(t)v\in 
\mathcal{H}\},
\end{equation*}%
and $A(t)v=A_{0}(t)v$ for all $v\in D(A(t)).$ Thus, (\ref{44})-(\ref{47})
can be written as the Cauchy problem 
\begin{eqnarray}
\frac{dY}{dt}(t)+A(t)Y(t) &=&f(t),\mbox{ a.e. }t\in (0,T),  \label{49} \\
Y(0) &=&Y_{0},  \notag
\end{eqnarray}%
and show further that it is well-posed. Since $A(t)$ is time dependent, the
existence proof relies on the result of Kato (see \cite{Kato}) extended by
Crandall and Pazy (see \cite{Crandall-Pazy-71}) for nonlinear evolution
equations. To this end we proceed to check, according e.g., to \cite%
{Crandall-Pazy-71}, the following properties of $A(t):$

(i) $\overline{D(A(t))}$ is independent of $t;$

(ii) $A(t)$ is quasi $m$-accretive on $\mathcal{H}$ for all $t\in \lbrack
0,T];$

(iii) For each $u\in \mathcal{H}$, $t\rightarrow J_{\lambda }(t)u$ is
Lipschitz from $[0,T]$ to $\mathcal{H}$, where $J_{\lambda }(t)$ is the
resolvent of $A(t).$

At point (i) we assert that $\overline{D(A(t))}=\mathcal{H}$ and this can
follow as a particular case of the proof given in \cite{cuiama-2},
Proposition 1, because one can note that $D(A(t))=\{v\in \mathcal{V};$ $%
v_{a}\in \mathcal{V}^{\prime },$ $v(0,x)=0,$ $v_{a}-\Delta v\in \mathcal{H}%
\}.$

Let $t$ be fixed. Let us compute 
\begin{eqnarray}
&&\left\langle A(t)v,v\right\rangle _{\mathcal{V}^{\prime },\mathcal{V}}=%
\frac{1}{2}\left\Vert v(a^{+})\right\Vert _{H}^{2}+\left\Vert \nabla
v\right\Vert _{\mathcal{H}}^{2}+\int_{0}^{a^{+}}\int_{O}(f_{1}(t)v^{2}+v%
\nabla v\cdot f_{2}(t))dxda  \notag \\
&&+\int_{0}^{a^{+}}\int_{\partial O}v^{2}f_{\Gamma }(t)d\sigma da\geq
\left\Vert v\right\Vert _{\mathcal{V}}^{2}-\left\Vert v\right\Vert _{%
\mathcal{H}}^{2}\left( \left\Vert f_{1}\right\Vert _{\infty }+\frac{1}{2}%
\left\Vert \nabla \cdot f_{2}\right\Vert _{\infty }+1\right) ,  \label{50}
\end{eqnarray}%
which shows that $A(t)$ is quasi accretive for $\lambda >\lambda
_{0}=\left\Vert f_{1}\right\Vert _{\infty }+\frac{1}{2}\left\Vert \nabla
\cdot f_{2}\right\Vert _{\infty }+1,$ where $\left\Vert f_{i}\right\Vert
_{\infty }=\left\Vert f_{i}\right\Vert _{L^{\infty }((0,T)\times
(0,a^{+})\times X_{i})},$ $X_{i}=O$, $\partial O,$ $i=1,2.$ Here, we used
the properties of $f_{1},$ $f_{2}$ and $f_{\Gamma }$ by (\ref{48}), and the
Gauss-Ostrogradski formula, namely 
\begin{eqnarray*}
&&\int_{0}^{a^{+}}\int_{O}v\nabla v\cdot f_{2}(t)dxda=\frac{1}{2}%
\int_{0}^{a^{+}}\int_{O}f_{2}(t)\cdot \nabla v^{2}dxda \\
&=&\frac{1}{2}\int_{0}^{a^{+}}\int_{O}(\nabla \cdot (f_{2}v^{2})-v^{2}\nabla
\cdot f_{2})dxda\leq \frac{1}{2}\left\Vert \nabla \cdot f_{2}\right\Vert
_{\infty }\left\Vert v\right\Vert _{\mathcal{H}}^{2}.
\end{eqnarray*}%
The operator is quasi $m$-accretive because the equation%
\begin{equation}
\lambda z+A(t)z=h  \label{51}
\end{equation}%
has a solution $z\in D(A(t)),$ for each $h\in \mathcal{H}$. Indeed, let us
introduce the linear Cauchy problem%
\begin{eqnarray}
\frac{dz}{da}(a)+B_{0}(t,a)z(a) &=&h(a)\mbox{, \ a.e. }a\in (0,a^{+}),
\label{52} \\
z(0) &=&0,  \notag
\end{eqnarray}%
where $B_{0}(t,a):V\rightarrow V^{\prime },$ 
\begin{eqnarray*}
&&\left\langle B_{0}(t,a)z,\psi \right\rangle _{V^{\prime },V} \\
&=&\int_{O}\left( \lambda z\psi +\nabla z\cdot \nabla \psi
+f_{1}(t,a,x)z\psi +\psi \nabla z\cdot f_{2}(t,a,x)\right) dx+\int_{\partial
O}zf_{\Gamma }(t,a,\sigma )\psi d\sigma ,\mbox{ }
\end{eqnarray*}%
for all $\psi \in V$ and $a\in \lbrack 0,a^{+}].$

Recall that $t$ is fixed. The operator $B_{0}(t,a)$ is bounded and $%
\left\langle B_{0}(t,a)z,\psi \right\rangle _{V^{\prime },V}\geq \left\Vert
z\right\Vert _{V}^{2}-\left\Vert z\right\Vert _{H}^{2}\left( \lambda
-\lambda _{0}\right) ,$ so that, by Lions' theorem (see \cite{lions}),
problem (\ref{52}) has a unique solution $z\in L^{2}(0,a^{+};V)\cap
W^{1,2}(0,a^{+};V^{\prime }).$ By (\ref{51}) $A(t)z=h-\lambda z\in \mathcal{H%
}$, hence $z\in D(A(t)).$

To prove (iii) we start from the resolvent equation (\ref{51}) which has a
unique solution, as seen before, denoted further by $z^{t}=(\lambda
I+A(t))^{-1}h$. Writing the difference between two equations (\ref{51})
considered for $A(t)$ and $A(s),$ 
\begin{equation}
\lambda (z^{t}-z^{s})+A(t)z^{t}-A(s)z^{s}=0,  \label{ec-dif}
\end{equation}%
setting $z:=z^{t}-z^{s}$ and multiplying scalarly in $\mathcal{H}$ by $z$ we
get 
\begin{eqnarray*}
&&\lambda \left\Vert z\right\Vert _{\mathcal{H}}^{2}+\left\Vert \nabla
z\right\Vert _{\mathcal{H}}^{2}+\frac{1}{2}\int_{O}\left\vert
z(a^{+})\right\vert ^{2}dx+\int_{0}^{a^{+}}\int_{\partial O}(f_{\Gamma
}(t)-f_{\Gamma }(s))z^{t}zd\sigma da+\int_{0}^{a^{+}}\int_{\partial
O}f_{\Gamma }(s)z^{2}d\sigma da \\
&&+\int_{0}^{a^{+}}\int_{O}(f_{1}(t)-f_{1}(s))zz^{t}dxda+\int_{0}^{a^{+}}%
\int_{O}f_{1}(s)z^{2}dxda \\
&&+\int_{0}^{a^{+}}\int_{O}(f_{2}(t)-f_{2}(s))z\cdot \nabla
z^{t}dxda+\int_{0}^{a^{+}}\int_{O}f_{2}(s)\cdot z\nabla zdxda=0.
\end{eqnarray*}%
By the regularity assumptions (\ref{48}) we have 
\begin{equation*}
\left\vert f_{i}(t)-f_{i}(s)\right\vert =\left\vert \int_{s}^{t}f_{i,\tau
}(\tau )d\tau \right\vert \leq \left\Vert f_{i,\tau }\right\Vert _{\infty
}\left\vert t-s\right\vert ,
\end{equation*}%
$f_{i,\tau }$ and $f_{\Gamma ,\tau }$ below being the partial derivatives of 
$f_{i},$ $i=1,2,$ and $f_{\Gamma },$ respectively, with respect to $t.$ Then,%
\begin{equation*}
\left\vert \int_{0}^{a^{+}}\int_{O}f_{2}(s)\cdot z\nabla zdxda\right\vert
\leq \frac{1}{2}\left\Vert \nabla \cdot f_{2}\right\Vert _{\infty
}\left\Vert z\right\Vert _{\mathcal{H}}^{2},
\end{equation*}%
\begin{eqnarray*}
&&\left\vert \int_{0}^{a^{+}}\int_{\partial O}\left\vert f_{\Gamma
}(t)-f_{\Gamma }(s)\right\vert z^{t}zd\sigma da\right\vert \leq \left\Vert
f_{\Gamma ,\tau }\right\Vert _{\infty }\left\vert t-s\right\vert
\int_{0}^{a^{+}}\left\Vert z(a)\right\Vert _{L^{2}(\partial O)}\left\Vert
z^{t}(a)\right\Vert _{L^{2}(\partial O)}da \\
&\leq &\left\Vert f_{\Gamma ,\tau }\right\Vert _{\infty }\left\vert
t-s\right\vert c_{tr}^{2}\int_{0}^{a^{+}}\left\Vert z(a)\right\Vert
_{V}\left\Vert z^{t}(a)\right\Vert _{V}da\leq \frac{1}{2}\left\Vert \nabla
z\right\Vert _{\mathcal{V}}^{2}+\frac{1}{2}\left\Vert f_{\Gamma ,\tau
}\right\Vert _{\infty }^{2}c_{tr}^{4}\left\vert t-s\right\vert
^{2}\left\Vert z^{t}\right\Vert _{\mathcal{V}}^{2},
\end{eqnarray*}%
where $c_{tr}$ is the constant in the trace theorem. Performing all
calculations we obtain 
\begin{eqnarray*}
&&\lambda \left\Vert z\right\Vert _{\mathcal{H}}^{2}+\frac{1}{2}\left\Vert
\nabla z\right\Vert _{\mathcal{V}}^{2}\leq (\lambda _{0}+2)\left\Vert
z\right\Vert _{\mathcal{H}}^{2} \\
&&+\left\vert t-s\right\vert ^{2}\left( \left\Vert f_{1,\tau }\right\Vert
_{\infty }^{2}\left\Vert z^{t}\right\Vert _{\mathcal{H}}^{2}+\left\Vert
f_{2,\tau }\right\Vert _{\infty }^{2}\left\Vert \nabla z^{t}\right\Vert _{%
\mathcal{H}}^{2}+\frac{1}{2}\left\Vert f_{\Gamma ,\tau }\right\Vert _{\infty
}^{2}c_{tr}^{4}\left\Vert z^{t}\right\Vert _{\mathcal{V}}^{2}\right) .
\end{eqnarray*}%
Relation%
\begin{equation*}
\left\langle A(t)v,v\right\rangle _{\mathcal{V}^{\prime },\mathcal{V}%
}=(A(t)v,v)_{\mathcal{H}}\geq \left\Vert v\right\Vert _{\mathcal{V}%
}^{2}-\lambda _{0}\left\Vert v\right\Vert _{\mathcal{H}}^{2}
\end{equation*}%
implies that 
\begin{equation*}
\left\Vert v\right\Vert _{\mathcal{V}}^{2}\leq (A(t)v,v)_{\mathcal{H}%
}+\lambda _{0}\left\Vert v\right\Vert _{\mathcal{H}}^{2}\leq \left\Vert
A(t)v\right\Vert _{\mathcal{H}}^{2}+(\lambda _{0}+1)\left\Vert v\right\Vert
_{\mathcal{H}}^{2}
\end{equation*}%
and so we deduce that 
\begin{equation*}
(\lambda -\lambda _{1})\left\Vert z\right\Vert _{\mathcal{H}}^{2}\leq
C\left\vert t-s\right\vert ^{2}(\left\Vert z^{t}\right\Vert _{\mathcal{H}%
}^{2}+\left\Vert A(t)z^{t}\right\Vert _{\mathcal{H}}^{2})
\end{equation*}%
with $C$ a constant and $\lambda >\lambda _{1}=\lambda _{0}+2.$

Let us note that $z^{t}=J_{\lambda }(t)h$, $z^{s}=J_{\lambda }(s)h$ and $%
\left\Vert A(t)z^{t}\right\Vert _{\mathcal{H}}=\left\Vert A_{\lambda
}(t)h\right\Vert _{\mathcal{H}}\leq \left\Vert A(t)h\right\Vert _{\mathcal{H}%
}.$ Thus, we obtain point (iii), as claimed.

Let $f\in W^{1,1}(0,T;\mathcal{H})$ and $y_{0}\in D(A(t))$. Then, the Cauchy
problem (\ref{49}) has a unique strong solution 
\begin{equation*}
Y\in C([0,T];\mathcal{H})\cap L^{\infty }(0,T;D(A(t)))\cap W^{1,\infty }(0,T;%
\mathcal{H)}.
\end{equation*}%
We multiply (\ref{49}) by $\psi \in W^{1,2}(0,T;\mathcal{H})\cap L^{2}(0,T;%
\mathcal{V})$, with $\psi _{a}\in L^{2}(0,T;\mathcal{V}^{\prime })$,
integrate over $(0,t)\times (0,a)\times O,$ and obtain 
\begin{eqnarray}
&&\int_{0}^{a}\int_{O}Y(t,s,x)\psi
(t,s,x)dxds-\int_{0}^{t}\int_{0}^{a}\int_{O}Y\psi _{\tau }dxdsd\tau
-\int_{0}^{a}\int_{O}Y_{0}\psi (0,s,x)dxds  \label{48-t} \\
&&+\int_{0}^{t}\int_{O}Y(\tau ,a,x)\psi (\tau ,a,x)dxd\tau
-\int_{0}^{t}\int_{0}^{a}\int_{O}Y\psi _{s}dxdsd\tau  \notag \\
&&+\int_{0}^{t}\int_{0}^{a}\int_{O}(\nabla Y\cdot \nabla \psi +f_{1}Y\psi
+\psi f_{2}\cdot \nabla Y)dxdsd\tau +\int_{0}^{t}\int_{0}^{a}\int_{\partial
O}f_{\Gamma }Y\psi d\sigma dsd\tau  \notag \\
&=&\int_{0}^{t}\int_{0}^{a}\int_{O}f\psi dxdsd\tau ,\mbox{ }  \notag
\end{eqnarray}%
which, in particular, for $t=T,$ $a=a^{+}$ and $\psi (T,a,x)=0$ yields (\ref%
{48-0}) with $E_{1}=E_{2}=f_{\Gamma }^{0}=0.$

Next, by setting $\psi =Y$ in (\ref{48-t}), we get, 
\begin{eqnarray*}
&&\left\Vert Y(t)\right\Vert _{\mathcal{H}}^{2}+\int_{0}^{t}\int_{O}Y^{2}(%
\tau ,a,x)dxd\tau +\int_{0}^{t}\left\Vert Y(\tau )\right\Vert _{\mathcal{V}%
}^{2}d\tau \\
&\leq &\left\Vert Y_{0}\right\Vert _{\mathcal{H}}^{2}+\int_{0}^{t}\left\Vert
f(\tau )\right\Vert _{\mathcal{V}^{\prime }}^{2}d\tau +2\left( 1+\left\Vert
f_{1}\right\Vert _{\infty }+\frac{1}{2}\left\Vert \nabla \cdot
f_{2}\right\Vert _{\infty }\right) \int_{0}^{t}\left\Vert Y(\tau
)\right\Vert _{\mathcal{H}}^{2}d\tau \\
&=&a(t)+a_{1}\int_{0}^{t}\left\Vert Y(\tau )\right\Vert _{\mathcal{H}%
}^{2}d\tau ,
\end{eqnarray*}%
with $a_{1}=2\left( 1+\left\Vert f_{1}\right\Vert _{\infty }+\frac{1}{2}%
\left\Vert \nabla \cdot f_{2}\right\Vert _{\infty }\right) $ and $%
a(t)=\left\Vert Y_{0}\right\Vert _{\mathcal{H}}^{2}+\int_{0}^{t}\left\Vert
f(\tau )\right\Vert _{\mathcal{V}^{\prime }}^{2}d\tau .$ By Gronwall's lemma
applied for $a(t)$ non-decreasing we get 
\begin{equation}
\left\Vert Y(t)\right\Vert _{\mathcal{H}}^{2}+\int_{0}^{t}\int_{O}Y^{2}(\tau
,a,x)dxd\tau +\int_{0}^{t}\left\Vert Y(\tau )\right\Vert _{\mathcal{V}%
}^{2}d\tau \leq e^{a_{1}T}\left( \left\Vert Y_{0}\right\Vert _{\mathcal{H}%
}^{2}+\int_{0}^{t}\left\Vert f(\tau )\right\Vert _{\mathcal{V}^{\prime
}}^{2}d\tau \right) .  \label{500}
\end{equation}%
Since the operator is linear we also have an estimate for the difference of
two solutions $Y_{1}$ and $Y_{2}$ corresponding to two pairs of data $%
\{Y_{0}^{1},f^{1}\}$ and $\{Y_{0}^{2},f^{2}\},$ 
\begin{eqnarray}
&&\left\Vert (Y_{1}-Y_{2})(t)\right\Vert _{\mathcal{H}}^{2}+\int_{0}^{t}%
\int_{O}(Y_{1}-Y_{2})^{2}(\tau ,a,x)dxd\tau +\int_{0}^{t}\left\Vert
(Y_{1}-Y_{2})(\tau )\right\Vert _{\mathcal{V}}^{2}d\tau  \label{501-1} \\
&\leq &\left( \left\Vert Y_{0}^{1}-Y_{0}^{2}\right\Vert _{\mathcal{H}%
}^{2}+\int_{0}^{t}\left\Vert (f^{1}-f^{2})(\tau )\right\Vert _{\mathcal{V}%
^{\prime }}^{2}d\tau \right) e^{a_{1}T}.  \notag
\end{eqnarray}

\noindent \textbf{Step 2.} Let $f\in L^{2}(0,T;\mathcal{H}),$ $Y_{0}\in 
\mathcal{H}$, $f_{\Gamma }^{0}\neq 0.$ Let us define $F_{\Gamma }(t)\in 
\mathcal{V}^{\prime }$ by 
\begin{equation}
\left\langle F_{\Gamma }(t),\psi \right\rangle _{\mathcal{V}^{\prime },%
\mathcal{V}}=-\int_{0}^{a^{+}}\int_{\partial O}f_{\Gamma }^{0}(t,s,\sigma
)\psi (a,\sigma )d\sigma da,\mbox{ for all }\psi \in \mathcal{V},\mbox{ a.e. 
}t\in (0,T),  \label{FGam}
\end{equation}%
and note that 
\begin{equation*}
\left\Vert F_{\Gamma }(t)\right\Vert _{\mathcal{V}^{\prime }}=\sup_{\psi \in 
\mathcal{V},\left\Vert \psi \right\Vert _{\mathcal{V}}\leq 1}\left\langle
F_{\Gamma }(t),\psi \right\rangle _{\mathcal{V}^{\prime },\mathcal{V}}\leq
c_{tr}\left\Vert f_{\Gamma }^{0}(t)\right\Vert
_{L^{2}(0,a^{+};L^{2}(\partial O))}.
\end{equation*}

Let $F_{\Gamma }^{n}\in W^{1,1}(0,T;\mathcal{H})$, $f^{n}\in W^{1,1}(0,T;%
\mathcal{H}),$ $Y_{0}^{n}\in D(A(t)),$ such that $F_{\Gamma }^{n}\rightarrow
F_{\Gamma }$ strongly in $L^{2}(0,T;\mathcal{V}^{\prime }),$ $%
f^{n}\rightarrow f$ strongly in $L^{2}(0,T;\mathcal{H}),$ and $%
Y_{0}^{n}\rightarrow Y_{0}$ strongly in $\mathcal{H},$ as $n\rightarrow
\infty .$ Thus, as $n\rightarrow \infty ,$ 
\begin{equation*}
\int_{0}^{t}\left\langle F_{\Gamma }^{n}(\tau ),\psi \right\rangle _{%
\mathcal{V}^{\prime },\mathcal{V}}d\tau \rightarrow \int_{0}^{t}\left\langle
F_{\Gamma }(\tau ),\psi \right\rangle _{\mathcal{V}^{\prime },\mathcal{V}%
}d\tau =-\int_{0}^{a^{+}}\int_{\partial O}f_{\Gamma }^{0}(t,s,\sigma )\psi
(a,\sigma )d\sigma da,\mbox{ for all }\psi \in \mathcal{V},
\end{equation*}%
and%
\begin{equation*}
\int_{0}^{t}\left\Vert F_{\Gamma }^{n}(\tau )\right\Vert _{\mathcal{V}%
^{\prime }}^{2}d\tau \rightarrow \int_{0}^{t}\left\Vert F_{\Gamma }(\tau
)\right\Vert _{\mathcal{V}^{\prime }}^{2}d\tau \leq
c_{tr}^{2}\int_{0}^{t}\left\Vert f_{\Gamma }^{0}(\tau )\right\Vert
_{L^{2}(0,a^{+};L^{2}(\partial O))}^{2}d\tau .
\end{equation*}

Let us consider problem (\ref{49}) with $f$ replaced by $f^{n}+F_{\Gamma
}^{n}.$ This has a unique solution $Y^{n}$ satisfying (\ref{48-t}), that is 
\begin{eqnarray}
&&\int_{0}^{a}\int_{O}Y^{n}(t,s,x)\psi
(t,s,x)dxds-\int_{0}^{t}\int_{0}^{a}\int_{O}Y^{n}\psi _{\tau }dxdsd\tau
-\int_{0}^{a}\int_{O}Y_{0}^{n}\psi (0,s,x)dxds  \label{50-0} \\
&&+\int_{0}^{t}\int_{O}Y^{n}(\tau ,a,x)\psi (\tau ,a,x)dxd\tau
-\int_{0}^{t}\int_{0}^{a}\int_{O}Y^{n}\psi _{s}dxdsd\tau   \notag \\
&&+\int_{0}^{t}\int_{0}^{a}\int_{O}(\nabla Y^{n}\cdot \nabla \psi
+f_{1}Y^{n}\psi +\psi f_{2}\cdot \nabla Y^{n})dxdsd\tau
+\int_{0}^{t}\int_{0}^{a}\int_{\partial O}f_{\Gamma }Y^{n}\psi d\sigma
dsd\tau   \notag \\
&=&\int_{0}^{t}\int_{0}^{a}\int_{O}(f^{n}+F_{\Gamma }^{n})\psi dxdsd\tau
=\int_{0}^{t}\int_{0}^{a}\int_{O}f^{n}\psi dxdsd\tau
+\int_{0}^{t}\int_{0}^{a}\left\langle F_{\Gamma }^{n}(\tau ,s),\psi (\tau
,s)\right\rangle _{V^{\prime },V}dsd\tau .  \notag
\end{eqnarray}%
Moreover, the solution satisfies estimate (\ref{500}), with $f$ replaced by $%
f^{n}+F_{\Gamma }^{n},$ 
\begin{eqnarray}
&&\left\Vert Y^{n}(t)\right\Vert _{\mathcal{H}}^{2}+\int_{0}^{t}%
\int_{O}(Y^{n})^{2}(\tau ,a,x)dxd\tau +\int_{0}^{t}\left\Vert Y^{n}(\tau
)\right\Vert _{\mathcal{V}}^{2}d\tau   \label{500-20} \\
&\leq &e^{a_{1}t}\left( \left\Vert Y_{0}^{n}\right\Vert _{\mathcal{H}%
}^{2}+2\int_{0}^{t}\left\Vert f^{n}(\tau )\right\Vert _{\mathcal{V}^{\prime
}}^{2}d\tau +2c_{tr}^{2}\int_{0}^{t}\left\Vert f_{\Gamma }^{0}(\tau
)\right\Vert _{L^{2}(0,a^{+};L^{2}(\partial O))}^{2}d\tau \right) ,  \notag
\end{eqnarray}%
and (\ref{501-1}), for the difference $Y_{1}^{n}-Y_{2}^{n},$ corresponding
to two sets of data, $\{Y_{0}^{i},f^{i},F_{\Gamma }^{i}\}_{i=1,2},$%
\begin{eqnarray}
&&\left\Vert (Y_{1}^{n}-Y_{2}^{n})(t)\right\Vert _{\mathcal{H}%
}^{2}+\int_{0}^{t}\int_{O}(Y_{1}^{n}-Y_{2}^{n})^{2}(\tau ,a,x)dxd\tau
+\int_{0}^{t}\left\Vert (Y_{1}^{n}-Y_{2}^{n})(\tau )\right\Vert _{\mathcal{V}%
}^{2}d\tau   \label{500-21} \\
&\leq &2e^{a_{1}T}\left( \left\Vert Y_{0}^{1n}-Y_{0}^{2n}\right\Vert _{%
\mathcal{H}}^{2}+\int_{0}^{t}\left\Vert (f^{1n}-f^{2n})(\tau )\right\Vert _{%
\mathcal{V}^{\prime }}^{2}d\tau +c_{tr}^{2}\int_{0}^{t}\left\Vert (F_{\Gamma
}^{1n}-F_{\Gamma }^{2n})(\tau )\right\Vert _{\mathcal{V}^{\prime }}^{2}d\tau
\right) .  \notag
\end{eqnarray}%
This particularized for $Y^{n}-Y^{m}$ gives 
\begin{eqnarray}
&&\left\Vert (Y^{n}-Y^{m})(t)\right\Vert _{\mathcal{H}}^{2}+\int_{0}^{t}%
\int_{O}(Y^{n}-Y^{m})^{2}(\tau ,a,x)dxd\tau +\int_{0}^{t}\left\Vert
(Y^{n}-Y^{m})(\tau )\right\Vert _{\mathcal{V}}^{2}d\tau   \label{500-2} \\
&\leq &2e^{a_{1}T}\left( \left\Vert Y_{0}^{n}-Y_{0}^{m}\right\Vert _{%
\mathcal{H}}^{2}+\int_{0}^{t}\left\Vert (f^{n}-f^{m})(\tau )\right\Vert _{%
\mathcal{V}^{\prime }}^{2}d\tau +c_{tr}^{2}\int_{0}^{t}\left\Vert (F_{\Gamma
}^{n}-F_{\Gamma }^{m})(\tau )\right\Vert _{\mathcal{V}^{\prime }}^{2}d\tau
\right)   \notag
\end{eqnarray}%
whence it follows that $\{Y^{n}\}_{n}$ is a Cauchy sequence in the spaces
indicated in (\ref{48-00}), therefore tending strongly to $Y$ in these
spaces. Moreover, by passing to the limit in (\ref{50-0}) we get that the
solution satisfies (\ref{48-t}) with the right-hand side 
\begin{equation*}
\int_{0}^{t}\int_{0}^{a}\int_{O}f\psi dxdsd\tau
-\int_{0}^{t}\int_{0}^{a}\int_{O}f_{\Gamma }^{0}\psi d\sigma dsd\tau .
\end{equation*}%
Next, (\ref{500-20}) and (\ref{500-21}) are preserved at limit, and imply 
\begin{eqnarray}
&&\left( \left\Vert Y_{0}\right\Vert _{\mathcal{H}}^{2}+\int_{0}^{t}\left%
\Vert f(\tau )+F_{\Gamma }(\tau )\right\Vert _{\mathcal{V}^{\prime
}}^{2}d\tau \right) e^{a_{1}T}  \label{500-22} \\
&\leq &c_{0}e^{a_{1}T}\left( \left\Vert Y_{0}\right\Vert _{\mathcal{H}%
}^{2}+\int_{0}^{t}\left\Vert f(\tau )\right\Vert _{\mathcal{H}}^{2}d\tau
+\int_{0}^{t}\left\Vert f_{\Gamma }^{0}(\tau )\right\Vert
_{L^{2}(0,a^{+};L^{2}(\partial O))}^{2}d\tau \right) ,  \notag
\end{eqnarray}%
(because $\left\Vert f(\tau )\right\Vert _{\mathcal{V}^{\prime }}\leq
\left\Vert f(\tau )\right\Vert _{\mathcal{H}})$ and 
\begin{eqnarray}
&&\left\Vert (Y-\overline{Y})(t)\right\Vert _{\mathcal{H}}^{2}+\int_{0}^{t}%
\int_{O}(Y-\overline{Y})^{2}(\tau ,a,x)d\tau dx+\int_{0}^{t}\left\Vert (Y-%
\overline{Y})(\tau )\right\Vert _{\mathcal{V}}^{2}d\tau   \label{508} \\
&\leq &c_{0}e^{a_{1}T}\left( \left\Vert Y_{0}^{1}-Y_{0}^{2}\right\Vert _{%
\mathcal{H}}^{2}+\int_{0}^{t}\left\Vert (f_{\Gamma }^{01}-f_{\Gamma
}^{02})(\tau )\right\Vert _{L^{2}(0,a^{+};L^{2}(\partial O))}^{2}d\tau
+\int_{0}^{t}\left\Vert (f^{1}-f^{2})(\tau )\right\Vert _{\mathcal{H}%
}^{2}d\tau \right) .  \notag
\end{eqnarray}%
The uniqueness is obvious. Here, $c_{0}$ is a constant depending on the
domain and dimension (via $c_{tr}).$

\noindent \textbf{Step 3.} Let $f\in L^{2}(0,T;\mathcal{H}),$ $Y_{0}\in 
\mathcal{H}$, $f_{\Gamma }^{0}\neq 0$ and let us consider the boundary
condition 
\begin{equation}
Y(t,0,x)=F(t,x)\mbox{ with\ }F\in L^{2}(0,T;H).  \label{58-0}
\end{equation}%
In a similar way as done at Step 2, we regularize all functions $f,Y_{0}$
and $F,$ for the last one choosing a sequence $F^{n}\in C^{2}([0,T]\times 
\overline{O}),$ such that $F^{n}(0,x)=0$ and $F^{n}\rightarrow F$ strongly
in $L^{2}(0,T;H).$ We have%
\begin{equation}
Y_{t}^{n}+Y_{a}^{n}-\Delta Y^{n}+f_{1}(t,a,x)Y^{n}+f_{2}(t,a,x)\cdot \nabla
Y^{n}=f^{n},\mbox{ in }(0,T)\times (0,a^{+})\times O,  \label{44-0}
\end{equation}%
\begin{equation}
-\nabla Y^{n}\cdot \nu =Y^{n}f_{\Gamma }+f_{\Gamma }^{0},\mbox{ in }%
(0,T)\times (0,a^{+})\times \partial O,  \label{45-0}
\end{equation}%
\begin{equation}
Y^{n}(t,0,x)=F^{n}(t,x),\mbox{ in }(0,T)\times O,  \label{46-0}
\end{equation}%
\begin{equation}
Y^{n}(0,a,x)=Y_{0}^{n}(a,x),\mbox{ in }(0,a^{+})\times O.  \label{47-0}
\end{equation}%
Homogenizing the boundary condition, by setting $Z:=Y-F^{n}$ we get the
system 
\begin{equation*}
Z_{t}+Z_{a}-\Delta Z+f_{1}(t,a,x)Z+f_{2}(t,a,x)\cdot \nabla Z=\widetilde{f}%
^{n},\mbox{ in }(0,T)\times (0,a^{+})\times O,
\end{equation*}%
\begin{equation*}
-\nabla Z\cdot \nu =Zf_{\Gamma }+f_{\Gamma }^{0}+\nabla F^{n}\cdot \nu
+f_{\Gamma }F^{n},\mbox{ in }(0,T)\times (0,a^{+})\times \partial O,
\end{equation*}%
\begin{equation*}
Z(t,0,x)=0,\mbox{ in }(0,T)\times O,
\end{equation*}%
\begin{equation*}
Z(0,a,x)=Z_{0}^{n}(a,x),\mbox{ in }(0,a^{+})\times O,
\end{equation*}%
where 
\begin{eqnarray*}
\widetilde{f}^{n} &=&f^{n}-F_{t}^{n}-F_{a}^{n}+\Delta
F^{n}-f_{1}F^{n}-f_{2}\cdot \nabla F^{n}\in W^{1,1}(0,T;\mathcal{H}), \\
\widetilde{f_{\Gamma }^{0}} &=&f_{\Gamma }^{0}+\nabla F^{n}\cdot \nu
+f_{\Gamma }F^{n}\in C^{1}([0,T];L^{2}((0,a^{+})\times \partial O)),\mbox{ }
\\
Z_{0}^{n} &=&Y_{0}^{n}-F^{n}(0,x)\in D(A(t)).
\end{eqnarray*}%
Denoting $\widetilde{F_{\Gamma }}^{n}=F_{\Gamma }^{n}+G^{n},$ where $%
F_{\Gamma }^{n}$ is the regularization of $F_{\Gamma }$ given by (\ref{FGam}%
) and 
\begin{equation*}
\left\langle G^{n}(t),\psi \right\rangle _{\mathcal{V}^{\prime },\mathcal{V}%
}=-\int_{0}^{a^{+}}\int_{\partial O}(\nabla F^{n}(t,a,\sigma )\cdot \nu
+f_{\Gamma }(t,a,\sigma )F^{n}(t,a,\sigma ))\psi (a,\sigma )d\sigma da,\mbox{
}
\end{equation*}%
we can write the Cauchy problem 
\begin{eqnarray*}
\frac{dZ}{dt}(t)+A(t)Z(t) &=&\widetilde{f}^{n}(t)+\widetilde{F_{\Gamma }}%
^{n}(t),\mbox{ a.e. }t\in (0,T), \\
Z(0) &=&Z_{0}^{n}.
\end{eqnarray*}%
Thus, we can apply Step 2 and assert that this new system has a unique
solution $Z^{n},$ satisfying 
\begin{eqnarray*}
&&\int_{0}^{a}\int_{O}Z^{n}(t,s,x)\psi
(t,s,x)dxds-\int_{0}^{t}\int_{0}^{a}\int_{O}Z^{n}\psi _{\tau }dxdsd\tau
-\int_{0}^{a}\int_{O}Z_{0}^{n}\psi (0,s,x)dxds \\
&&+\int_{0}^{t}\int_{O}Z^{n}(\tau ,a,x)\psi (\tau ,a,x)dxd\tau
-\int_{0}^{t}\int_{0}^{a}\int_{O}Z^{n}\psi _{s}dxdsd\tau \\
&&+\int_{0}^{t}\int_{0}^{a}\int_{O}(\nabla Z^{n}\cdot \nabla \psi
+f_{1}Z^{n}+f_{2}\psi \cdot \nabla Z^{n})dxdsd\tau
+\int_{0}^{t}\int_{0}^{a}\int_{\partial O}f_{\Gamma }Z^{n}\psi d\sigma
dsd\tau \\
&=&\int_{0}^{t}\int_{0}^{a}\int_{O}\left( \widetilde{f}^{n}+\widetilde{%
F_{\Gamma }}^{n}\right) \psi dxdsd\tau .
\end{eqnarray*}%
Making some computations for going back to $Y^{n}=Z^{n}+F^{n}$ we get that
it satisfies 
\begin{eqnarray}
&&\int_{0}^{a}\int_{O}Y^{n}(t,s,x)\psi
(t,s,x)dxds-\int_{0}^{t}\int_{0}^{a}\int_{O}Y^{n}\psi _{\tau }dxdsd\tau
-\int_{0}^{a}\int_{O}Y_{0}^{n}\psi (0,s,x)dxds  \label{505-0} \\
&&+\int_{0}^{t}\int_{O}Y^{n}(\tau ,a,x)\psi (\tau ,a,x)dxd\tau
-\int_{0}^{t}\int_{0}^{a}\int_{O}Y^{n}\psi _{s}dxdsd\tau  \notag \\
&&-\int_{0}^{t}\int_{O}F^{n}(\tau ,x)\psi (\tau ,0,x)dxd\tau
+\int_{0}^{t}\int_{0}^{a}\int_{O}(\nabla Y^{n}\cdot \nabla \psi
+f_{1}Y^{n}+f_{2}\psi \cdot \nabla Y^{n})dxdsd\tau +  \notag \\
&&\int_{0}^{t}\int_{0}^{a}\int_{\partial O}(f_{\Gamma }Y^{n}\psi +f_{\Gamma
}^{0}\psi )d\sigma dsd\tau =\int_{0}^{t}\int_{0}^{a}\int_{O}f^{n}\psi
dxdsd\tau +\int_{0}^{t}\left\langle F_{\Gamma }^{n}(\tau ),\psi
\right\rangle _{\mathcal{V}^{\prime },\mathcal{V}}d\tau .  \notag
\end{eqnarray}%
Therefore, we obtain the estimates 
\begin{eqnarray}
&&\left\Vert Y^{n}(t)\right\Vert _{\mathcal{H}}^{2}+\int_{0}^{t}%
\int_{O}(Y^{n})^{2}(\tau ,a,x)dxd\tau +\int_{0}^{t}\left\Vert Y^{n}(\tau
)\right\Vert _{\mathcal{V}}^{2}d\tau  \label{505} \\
&\leq &c_{0}e^{a_{1}T}\left( \left\Vert Y_{0}^{n}\right\Vert _{\mathcal{H}%
}^{2}+\int_{0}^{t}\int_{O}(F^{n})^{2}dxd\tau +\int_{0}^{t}\left\Vert
f^{n}(\tau )\right\Vert _{\mathcal{V}^{\prime }}^{2}d\tau
+\int_{0}^{t}\left\Vert F_{\Gamma }^{n}(\tau )\right\Vert _{\mathcal{V}%
^{\prime }}^{2}d\tau \right)  \notag
\end{eqnarray}%
and 
\begin{eqnarray}
&&\left\Vert (Y^{n}-\overline{Y^{n}})(t)\right\Vert _{\mathcal{H}%
}^{2}+\int_{0}^{t}\int_{O}(Y^{n}-\overline{Y^{n}})^{2}(\tau ,a,x)dxd\tau
+\int_{0}^{t}\left\Vert (Y^{n}-\overline{Y^{n}})(\tau )\right\Vert _{%
\mathcal{V}}^{2}d\tau  \label{506} \\
&\leq &c_{0}e^{a_{1}T}\left( \left\Vert Y_{0}^{n}-\overline{Y_{0}^{n}}%
\right\Vert _{\mathcal{H}}^{2}+\int_{0}^{t}\int_{O}(F^{n}-\overline{F^{n}}%
)^{2}dxd\tau \right.  \notag \\
&&\left. +\int_{0}^{t}\left\Vert (f^{n}-\overline{f^{n}})(\tau )\right\Vert
_{\mathcal{H}}^{2}d\tau +\int_{0}^{t}\left\Vert (F_{\Gamma }^{n}-\overline{%
F_{\Gamma }^{n}})(\tau )\right\Vert _{\mathcal{V}^{\prime }}d\tau \right) , 
\notag
\end{eqnarray}%
where $\overline{Y^{n}}$ is the solution corresponding to $\{\overline{Y_{0}}%
,$ $f_{\Gamma },$ $\overline{f_{\Gamma }^{0}},$ $\overline{F},$ $\overline{f}%
\}$ and $a_{1}$ depends on the problem parameters ($\left\Vert
f_{1}\right\Vert _{\infty },$ $\left\Vert \nabla \cdot f_{2}\right\Vert
_{\infty })$ and $T$. Here, $f_{\Gamma }$ is the same for both solutions.

Arguing as before, we get that $\{Y^{n}\}_{n}$ is Cauchy in $C([0,T];%
\mathcal{H})\cap C([0,a^{+}];L^{2}(0,T;H))\cap L^{2}(0,T;\mathcal{V}),$
hence $Y^{n}\rightarrow Y$ strongly in these spaces as $n\rightarrow \infty
, $ so that by passing to the limit in (\ref{505-0}) we obtain 
\begin{eqnarray}
&&\int_{0}^{a}\int_{O}Y(t,s,x)\psi
(t,s,x)dxda-\int_{0}^{t}\int_{0}^{a}\int_{O}Y\psi _{\tau }dxdsd\tau
-\int_{0}^{a}\int_{O}Y_{0}\psi (0,s,x)dxds  \label{502} \\
&&+\int_{0}^{t}\int_{O}Y(\tau ,a,x)\psi (\tau ,a,x)dxd\tau
-\int_{0}^{t}\int_{0}^{a}\int_{O}Y\psi _{s}dxdsd\tau  \notag \\
&&-\int_{0}^{t}\int_{O}F(\tau ,x)\psi (\tau ,0,x)dxd\tau
+\int_{0}^{t}\int_{0}^{a}\int_{O}(\nabla Y\cdot \nabla \psi
+f_{1}Y+f_{2}\psi \cdot \nabla Y)dxdsd\tau +  \notag \\
&&\int_{0}^{t}\int_{0}^{a}\int_{\partial O}(f_{\Gamma }Y\psi +f_{\Gamma
}^{0}\psi )d\sigma dsd\tau =\int_{0}^{t}\int_{0}^{a}\int_{O}f\psi dxdsd\tau .
\notag
\end{eqnarray}%
Setting $t=T$ and $a=a^{+}$ and taking $\psi (T,a,x)=0,$ we get that system (%
\ref{44}), (\ref{45}), (\ref{47}), (\ref{58-0}) has a solution.\ By passing
to the limit as $n\rightarrow \infty ,$ in estimates (\ref{505}) and (\ref%
{506}) we get 
\begin{eqnarray}
&&\left\Vert Y(t)\right\Vert _{\mathcal{H}}^{2}+\int_{0}^{t}\int_{O}Y^{2}(%
\tau ,a,x)dxd\tau +\int_{0}^{t}\left\Vert Y(\tau )\right\Vert _{\mathcal{V}%
}^{2}d\tau  \label{505-00} \\
&\leq &c_{0}e^{a_{1}T}\left( \left\Vert Y_{0}\right\Vert _{\mathcal{H}%
}^{2}+\int_{0}^{t}\int_{O}F^{2}dxd\tau +\int_{0}^{t}\left\Vert f(\tau
)\right\Vert _{\mathcal{H}}^{2}d\tau +\int_{0}^{t}\left\Vert f_{\Gamma
}^{0}(\tau )\right\Vert _{L^{2}(0,a^{+};L^{2}(\partial O))}^{2}d\tau \right)
:=\overline{C}  \notag
\end{eqnarray}%
and 
\begin{eqnarray}
&&\left\Vert (Y-\overline{Y})(t)\right\Vert _{\mathcal{H}}^{2}+\int_{0}^{t}%
\int_{O}(Y-\overline{Y})^{2}(\tau ,a,x)dxd\tau +\int_{0}^{t}\left\Vert (Y-%
\overline{Y})(\tau )\right\Vert _{\mathcal{V}}^{2}d\tau  \label{506-00} \\
&\leq &c_{0}e^{a_{1}T}\left( \left\Vert Y_{0}-\overline{Y_{0}}\right\Vert _{%
\mathcal{H}}^{2}+\int_{0}^{t}\int_{O}(F-\overline{F})^{2}dxd\tau \right. 
\notag \\
&&\left. +\int_{0}^{t}\left\Vert (f-\overline{f})(\tau )\right\Vert _{%
\mathcal{V}^{\prime }}^{2}d\tau +\int_{0}^{t}\left\Vert (f_{\Gamma }^{0}-%
\overline{f_{\Gamma }^{0}})(\tau )\right\Vert _{L^{2}(0,a^{+};L^{2}(\partial
O))}^{2}d\tau \right) ,  \notag
\end{eqnarray}%
corresponding to two sets of data $\{Y_{0},f,f_{\Gamma }^{0},F\}$ and $\{%
\overline{Y_{0}},\overline{f},\overline{f_{\Gamma }^{0}},\overline{F}\}$ and
to the same $f_{1},$ $f_{2}$ and $f_{\Gamma }$. Again (\ref{506-00}) ensures
the uniqueness.

For a later use we deduce the estimate for the difference of two solutions $%
Y_{1}$ and $Y_{2}$ corresponding to two completely different sets of data $%
\{Y_{0}^{1},f_{1}^{1},f_{2}^{1},f_{\Gamma }^{1},f^{1},F^{1}\}$ and $%
\{Y_{0}^{2},f_{1}^{2},f_{2}^{2},f_{\Gamma }^{2},f^{2},F^{2}\},$ computing
first the estimate for the regular solutions and then passing to the limit.
For simplicity we do not indicate the superscript $n$ for the regularized
solutions in the following computations. We have%
\begin{eqnarray}
&&\frac{1}{2}\left\Vert (Y_{1}-Y_{2})(t)\right\Vert _{\mathcal{H}%
}^{2}+\int_{0}^{t}\int_{O}(Y_{1}-Y_{2})^{2}(\tau ,a,x)dxd\tau
+\int_{0}^{t}\left\Vert (Y_{1}-Y_{2})(\tau )\right\Vert _{\mathcal{V}%
}^{2}d\tau  \label{dif-total} \\
&&+\int_{0}^{t}\int_{0}^{a}\int_{\partial O}(Y^{1}-Y^{2)^{2}}f_{\Gamma
}^{2}d\sigma dsd\tau  \notag \\
&\leq &\frac{1}{2}\left\Vert Y_{0}^{1}-Y_{0}^{2}\right\Vert _{\mathcal{H}%
}^{2}+6\int_{0}^{t}\left\Vert (f^{1}-f^{2})(\tau )\right\Vert _{\mathcal{V}%
^{\prime }}^{2}d\tau +\frac{1}{6}\int_{0}^{t}\left\Vert (Y_{1}-Y_{2})(\tau
)\right\Vert _{\mathcal{V}}^{2}d\tau  \notag \\
&&+\int_{0}^{t}\int_{0}^{a}\int_{O}\left( \left\vert
f_{1}^{1}-f_{1}^{2}\right\vert \left\vert Y^{1}\right\vert \left\vert
Y^{1}-Y^{2}\right\vert +\left\vert Y^{1}-Y^{2}\right\vert ^{2}\left\vert
f_{1}^{2}\right\vert \right) dxdsd\tau  \notag \\
&&+\int_{0}^{t}\int_{0}^{a}\int_{O}\left( \left\vert
f_{2}^{1}-f_{2}^{2}\right\vert \left\vert \nabla Y^{1}\right\vert \left\vert
Y^{1}-Y^{2}\right\vert +\frac{1}{2}\left\vert Y^{1}-Y^{2}\right\vert
^{2}\left\vert \nabla \cdot f_{2}^{2}\right\vert \right) dxdsd\tau  \notag \\
&&+\frac{1}{2}\int_{0}^{t}\int_{O}(F^{1}-F^{2})^{2}dxd\tau
+\int_{0}^{t}\left\Vert (Y_{1}-Y_{2})(\tau )\right\Vert _{\mathcal{H}%
}^{2}d\tau  \notag \\
&&+\int_{0}^{t}\int_{0}^{a}\int_{\partial O}\left\vert f_{\Gamma
}^{1}-f_{\Gamma }^{2}\right\vert \left\vert Y^{1}\right\vert \left\vert
Y^{1}-Y^{2}\right\vert d\sigma dsd\tau
+\int_{0}^{t}\int_{0}^{a}\int_{\partial O}\left\vert f_{\Gamma
}^{01}-f_{\Gamma }^{02}\right\vert \left\vert Y^{1}-Y^{2}\right\vert d\sigma
dsd\tau .  \notag
\end{eqnarray}%
Further we have%
\begin{eqnarray*}
&&\frac{1}{2}\left\Vert (Y_{1}-Y_{2})(t)\right\Vert _{\mathcal{H}}^{2}+\frac{%
1}{2}\int_{0}^{t}\int_{O}(Y_{1}-Y_{2})^{2}(\tau ,a,x)dxd\tau
+\int_{0}^{t}\left\Vert (Y_{1}-Y_{2})(\tau )\right\Vert _{\mathcal{V}%
}^{2}d\tau \\
&&+\int_{0}^{t}\int_{0}^{a}\int_{\partial O}\left\vert
Y^{1}-Y^{2}\right\vert ^{2}\left\vert f_{\Gamma }^{2}\right\vert d\sigma
dsd\tau \\
&\leq &\frac{1}{2}\left\Vert Y_{0}^{1}-Y_{0}^{2}\right\Vert _{\mathcal{H}%
}^{2}+6\int_{0}^{t}\left\Vert (f^{1}-f^{2})(\tau )\right\Vert _{\mathcal{H}%
}^{2}d\tau +\frac{1}{6}\int_{0}^{t}\left\Vert (Y_{1}-Y_{2})(\tau
)\right\Vert _{\mathcal{V}}^{2}d\tau \\
&&+\left\Vert f_{1}^{1}-f_{1}^{2}\right\Vert _{\infty }^{2}\left\Vert
Y^{1}(\tau )\right\Vert _{\mathcal{H}}^{2}+\int_{0}^{t}\left\Vert
(Y^{1}-Y^{2})(\tau )\right\Vert _{\mathcal{H}}^{2}d\tau +\left\Vert
f_{1}^{2}\right\Vert _{\infty }\int_{0}^{t}\left\Vert (Y^{1}-Y^{2})(\tau
)\right\Vert _{\mathcal{H}}^{2}d\tau \\
&&+\left\Vert f_{2}^{1}-f_{2}^{2}\right\Vert _{\infty }^{2}\left\Vert \nabla
Y^{1}(\tau )\right\Vert _{\mathcal{H}}^{2}+\int_{0}^{t}\left\Vert
(Y^{1}-Y^{2})(\tau )\right\Vert _{\mathcal{H}}^{2}d\tau +\frac{1}{2}%
\left\Vert \nabla \cdot f_{2}^{2}\right\Vert _{\infty
}\int_{0}^{t}\left\Vert (Y^{1}-Y^{2})(\tau )\right\Vert _{\mathcal{H}%
}^{2}d\tau \\
&&+\frac{1}{2}\int_{0}^{t}\int_{O}(F^{1}-F^{2})^{2}dxd\tau
+\int_{0}^{t}\left\Vert (Y^{1}-Y^{2})(\tau )\right\Vert _{\mathcal{H}%
}^{2}d\tau \\
&&+\frac{1}{6}\int_{0}^{t}\left\Vert (Y_{1}-Y_{2})(\tau )\right\Vert _{%
\mathcal{V}}^{2}d\tau +6c_{tr}^{4}\left\Vert f_{\Gamma }^{1}-f_{\Gamma
}^{2}\right\Vert _{\infty }^{2}\left\Vert Y^{1}(\tau )\right\Vert _{\mathcal{%
V}}^{2} \\
&&+\frac{1}{6}\int_{0}^{t}\left\Vert (Y_{1}-Y_{2})(\tau )\right\Vert _{%
\mathcal{V}}^{2}d\tau +6c_{tr}^{2}\int_{0}^{t}\left\Vert (f_{\Gamma
}^{01}-f_{\Gamma }^{02})(\tau )\right\Vert _{L^{2}(0,a^{+};L^{2}(\partial
O))}^{2}d\tau ,
\end{eqnarray*}%
which implies by (\ref{505-00}) 
\begin{eqnarray}
&&\left\Vert (Y_{1}-Y_{2})(t)\right\Vert _{\mathcal{H}}^{2}+\int_{0}^{t}%
\int_{O}(Y_{1}-Y_{2})^{2}(\tau ,a,x)dxd\tau +\int_{0}^{t}\left\Vert
(Y_{1}-Y_{2})(\tau )\right\Vert _{\mathcal{V}}^{2}d\tau  \label{dif-total-1}
\\
&\leq &c_{1}\left\{ \left\Vert Y_{0}^{1}-Y_{0}^{2}\right\Vert _{\mathcal{H}%
}^{2}+\int_{0}^{t}\left\Vert (f^{1}-f^{2})(\tau )\right\Vert _{\mathcal{H}%
}^{2}d\tau +\int_{0}^{t}\int_{O}(F^{1}-F^{2})^{2}dxd\tau \right\}  \notag \\
&&+c_{0}\int_{0}^{t}\left\Vert (f_{\Gamma }^{01}-f_{\Gamma }^{02})(\tau
)\right\Vert _{L^{2}(0,a^{+};L^{2}(\partial O))}^{2}d\tau +c_{0}\overline{C}%
_{1}^{2}\left( \left\Vert f_{1}^{1}-f_{1}^{2}\right\Vert _{\infty
}^{2}+\left\Vert f_{2}^{1}-f_{2}^{2}\right\Vert _{\infty }^{2}+\left\Vert
f_{\Gamma }^{1}-f_{\Gamma }^{2}\right\Vert _{\infty }^{2}\right)  \notag \\
&&+c_{1}\left( 1+\left\Vert f_{1}^{2}\right\Vert _{\infty }+\left\Vert
\nabla \cdot f_{2}^{2}\right\Vert _{\infty }\right) \int_{0}^{t}\left\Vert
(Y^{1}-Y^{2})(\tau )\right\Vert _{\mathcal{H}}^{2}d\tau ,  \notag
\end{eqnarray}%
where $\overline{C}_{1}$ is $\overline{C}$ given by (\ref{505-00})
corresponding to the functions indexed by $1$ and $c_{0}$ is a constant
depending on $c_{tr}.$

\noindent \textbf{Step 4.} Let us consider the complete system (\ref{44})-(%
\ref{47}). We shall apply the Banach fixed point theorem in the space $%
C([0,T;\mathcal{H}).$ Let us fix $\zeta \in C([0,T;\mathcal{H})$ and
consider the problem 
\begin{equation}
v_{t}+v_{a}-\Delta v+f_{1}(t,a,x)v+f_{2}(t,a,x)\cdot \nabla v=f^{\zeta
}(t,a,x),\mbox{ in }(0,T)\times (0,a^{+})\times O,  \label{40}
\end{equation}%
\begin{equation}
-\nabla v\cdot \nu =f_{\Gamma }(t,a,x)v+f_{\Gamma }^{0}(t,a,x),\mbox{ in }%
(0,T)\times (0,a^{+})\times \partial O,  \label{41}
\end{equation}%
\begin{equation}
v(t,0,x)=F^{\zeta }(t,x),\mbox{ in }(0,T)\times O,  \label{42}
\end{equation}%
\begin{equation}
v(0,a,x)=Y_{0}(a,x),\mbox{ in }(0,a^{+})\times O,  \label{43}
\end{equation}%
where%
\begin{equation*}
f^{\zeta }(t,a,x)=-E_{1}(t,a,x;\zeta ),\mbox{ \ }F^{\zeta
}(t,x)=\int_{0}^{a^{+}}E_{2}(t,a,x;\zeta )da.
\end{equation*}%
Note that $f^{\zeta }\in L^{2}(0,T;\mathcal{H})$ and $F^{\zeta }\in
L^{2}(0,T;H)$ and so we are entitled to apply Step 3 to find that system (%
\ref{40})-(\ref{43}) has a unique solution 
\begin{equation*}
v^{\zeta }\in C([0,T];\mathcal{H})\cap C([0,a^{+}];L^{2}(0,T;H))\cap
L^{2}(0,T;\mathcal{V})
\end{equation*}%
satisfying (\ref{502}) and estimate (\ref{505-00}), that is 
\begin{eqnarray}
&&\left\Vert v^{\zeta }(t)\right\Vert _{\mathcal{H}}^{2}+\left\Vert v^{\zeta
}(a)\right\Vert _{L^{2}(0,T;H)}^{2}+\int_{0}^{t}\left\Vert v^{\zeta
}(t)\right\Vert _{\mathcal{V}}^{2}d\tau  \label{43-0} \\
&\leq &c_{0}e^{a_{1}T}\left( \left\Vert y_{0}\right\Vert _{\mathcal{H}%
}^{2}+\int_{0}^{t}\int_{O}(F^{\zeta })^{2}(\tau ,x)dxd\tau
+\int_{0}^{t}\left\Vert f_{\Gamma }^{0}(\tau )\right\Vert
_{L^{2}(0,a^{+};L^{2}(\partial O))}^{2}d\tau +\int_{0}^{t}\left\Vert
E_{1}(\tau )\right\Vert _{\mathcal{H}}^{2}d\tau \right)  \notag \\
&\leq &c_{0}e^{a_{1}T}\left( \left\Vert y_{0}\right\Vert _{\mathcal{H}%
}^{2}+\left( a^{+}m_{\infty }^{2}+\mu _{\infty }^{2}\right)
\int_{0}^{t}\left\Vert \zeta (\tau )\right\Vert _{\mathcal{H}}^{2}d\tau
+\int_{0}^{T}\left\Vert f_{\Gamma }^{0}(t)\right\Vert
_{L^{2}(0,a^{+};L^{2}(\partial O))}^{2}dt\right) .  \notag
\end{eqnarray}%
For the passage to the last line in (\ref{43-0}) we used the properties (\ref%
{cE}) for $E_{1}$ and $E_{2},$ e.g., 
\begin{eqnarray*}
\int_{0}^{t}\int_{O}(F^{\zeta })^{2}(\tau ,x)dxd\tau
&=&\int_{0}^{t}\int_{O}\left( \int_{0}^{a^{+}}E_{2}(t,a,x;\zeta )da\right)
^{2}dxd\tau \leq
a^{+}\int_{0}^{t}\int_{O}\int_{0}^{a^{+}}E_{2}^{2}(t,a,x;\zeta )dadxd\tau \\
&\leq &a^{+}m_{\infty }^{2}\int_{0}^{t}\left\Vert \zeta (\tau )\right\Vert _{%
\mathcal{H}}^{2}d\tau .
\end{eqnarray*}%
Then, we define $\Psi :\mathcal{M}\rightarrow C([0,T];\mathcal{H})$ which
maps $\zeta \in \mathcal{M}$ into the solution $v^{\zeta }$ to (\ref{40})-(%
\ref{43}). Obviously, $\Psi (\mathcal{M})\subset \mathcal{M}$ and we show
that $\Psi $ is a contraction on $\mathcal{M}.$ Indeed, let $v^{\zeta }$ and 
$v^{\overline{\zeta }}$ be two solutions to (\ref{40})-(\ref{43})
corresponding to $\zeta $ and $\overline{\zeta }.$ Then, by (\ref{506-00}),
the estimate of the difference of these solutions reads%
\begin{eqnarray}
&&\left\Vert v^{\zeta }(t)-v^{\overline{\zeta }}(t)\right\Vert _{\mathcal{H}%
}^{2}+\left\Vert v^{\zeta }(a)-v^{\overline{\zeta }}(a)\right\Vert
_{L^{2}(0,T;H)}^{2}+\int_{0}^{t}\left\Vert v^{\zeta }(\tau )-v^{\overline{%
\zeta }}(\tau )\right\Vert _{\mathcal{V}}^{2}d\tau  \label{43-dif} \\
&\leq &c_{0}e^{a_{1}T}\left( \int_{0}^{t}\int_{O}(F^{\zeta }(\tau ,x)-F^{%
\overline{\zeta }}(\tau ,x))^{2}dxd\tau
+\int_{0}^{t}\int_{0}^{a^{+}}\int_{O}(E_{1}(\tau ,a,x;\zeta )-E_{1}(\tau
,a,x;\overline{\zeta }))^{2}dxdad\tau \right)  \notag \\
&\leq &c_{0}e^{a_{1}T}\left( a^{+}L_{2}^{2}\int_{0}^{t}\left\Vert (\zeta -%
\overline{\zeta })(\tau )\right\Vert _{\mathcal{H}}^{2}d\tau
+L_{1}^{2}\int_{0}^{t}\left\Vert (\zeta -\overline{\zeta })(\tau
)\right\Vert _{\mathcal{H}}^{2}d\tau \right) \leq C\int_{0}^{t}\left\Vert
(\zeta -\overline{\zeta })(\tau )\right\Vert _{\mathcal{H}}^{2}d\tau . 
\notag
\end{eqnarray}%
Considering now the norm $\left\Vert v^{\zeta }\right\Vert
_{B}=\sup\limits_{t\in \lbrack 0,T]}\left( e^{-\gamma _{0}t}\left\Vert
v^{\zeta }(t)\right\Vert _{\mathcal{H}}\right) $ which is equivalent with
the standard norm in $C([0,T];\mathcal{H}),$ it follows by some calculations
that 
\begin{eqnarray*}
\left\Vert v^{\zeta }-v^{\overline{\zeta }}\right\Vert _{B}^{2} &\leq
&Ce^{-2\gamma _{0}t}\int_{0}^{t}e^{2\gamma _{0}s}\left\Vert \zeta -\overline{%
\zeta }\right\Vert _{B}^{2}ds\leq \frac{C}{2\gamma _{0}}(1-e^{-2\gamma
_{0}t})\left\Vert \zeta -\overline{\zeta }\right\Vert _{B}^{2} \\
&\leq &\frac{C}{2\gamma _{0}}\left\Vert \zeta -\overline{\zeta }\right\Vert
_{B}^{2}
\end{eqnarray*}%
which proves, by a suitable choice $2\gamma _{0}>C,$ that $\Psi $ is a
contraction on $\mathcal{M}$. Then, $\Psi $ has a fixed point, $\Psi (\zeta
)=\zeta =v^{\zeta },$ which is the unique solution to (\ref{40})-(\ref{43}).
Thus, $v^{\zeta }$ turns out to solve (\ref{44})-(\ref{47}) and actually it
can be denoted by $Y.$

Finally, assuming that on the right-hand side of (\ref{44}) we have $%
f-E_{1}(t,a,x;Y),$ we get by using (\ref{dif-total-1}) that%
\begin{eqnarray}
&&\left\Vert (Y_{1}-Y_{2})(t)\right\Vert _{\mathcal{H}}^{2}+\int_{0}^{t}%
\int_{O}(Y_{1}-Y_{2})^{2}(\tau ,a,x)dxd\tau +\int_{0}^{t}\left\Vert
(Y_{1}-Y_{2})(\tau )\right\Vert _{\mathcal{V}}^{2}d\tau   \label{600} \\
&\leq &c_{1}\left\{ \left\Vert Y_{0}^{1}-Y_{0}^{2}\right\Vert _{\mathcal{H}%
}^{2}+\int_{0}^{t}\left\Vert (f^{1}-f^{2})(\tau )\right\Vert _{\mathcal{H}%
}^{2}d\tau \right.   \notag \\
&&+\int_{0}^{t}\int_{0}^{a^{+}}\int_{O}(E_{1}(\tau ,a,x;Y_{1})-E_{1}(\tau
,a,x;Y_{2}))^{2}dxdad\tau   \notag \\
&&\left. +\int_{0}^{t}\int_{O}\left( \int_{0}^{a^{+}}(E_{2}(\tau
,a,x;Y_{1})-E_{2}(\tau ,a,x;Y_{2}))da\right) ^{2}dxd\tau \right\}   \notag \\
&&+c_{0}\overline{C}_{1}^{2}\left( \left\Vert f_{1}^{1}-f_{1}^{2}\right\Vert
_{\infty }^{2}+\left\Vert f_{2}^{1}-f_{2}^{2}\right\Vert _{\infty
}^{2}+\left\Vert f_{\Gamma }^{1}-f_{\Gamma }^{2}\right\Vert _{\infty
}^{2}\right)   \notag \\
&&+c_{0}\int_{0}^{t}\left\Vert (f_{\Gamma }^{01}-f_{\Gamma }^{02})(\tau
)\right\Vert _{L^{2}(0,a^{+};L^{2}(\partial O))}^{2}d\tau   \notag \\
&&+c_{1}\left( 1+\left\Vert f_{1}^{2}\right\Vert _{\infty }+\left\Vert
f_{2}^{2}\right\Vert _{\infty }\right) \int_{0}^{t}\left\Vert
(Y^{1}-Y^{2})(\tau )\right\Vert _{\mathcal{H}}^{2}d\tau ,  \notag
\end{eqnarray}%
whence using the Lipschitz property of $E_{i}$, with $L_{1}$ and $L_{2}$
given by (\ref{Li}), we get%
\begin{eqnarray*}
&&\left\Vert (Y_{1}-Y_{2})(t)\right\Vert _{\mathcal{H}}^{2}+\int_{0}^{t}%
\int_{O}(Y_{1}-Y_{2})^{2}(\tau ,a,x)dxd\tau +\int_{0}^{t}\left\Vert
(Y_{1}-Y_{2})(\tau )\right\Vert _{\mathcal{V}}^{2}d\tau  \\
&\leq &c_{0}\left\{ \left\Vert Y_{0}^{1}-Y_{0}^{2}\right\Vert _{\mathcal{H}%
}^{2}+\int_{0}^{t}\left\Vert (f^{1}-f^{2})(\tau )\right\Vert _{\mathcal{H}%
}^{2}d\tau +\int_{0}^{t}\left\Vert (f_{\Gamma }^{01}-f_{\Gamma }^{02})(\tau
)\right\Vert _{L^{2}(0,a^{+};L^{2}(\partial O))}^{2}d\tau \right.  \\
&&\left. +\overline{C}_{1}^{2}\left( \left\Vert
f_{1}^{1}-f_{1}^{2}\right\Vert _{\infty }^{2}+\left\Vert
f_{2}^{1}-f_{2}^{2}\right\Vert _{\infty }^{2}+\left\Vert f_{\Gamma
}^{1}-f_{\Gamma }^{2}\right\Vert _{\infty }^{2}\right) \right\}  \\
&&+c_{1}\left( 1+\left\Vert f_{1}^{2}\right\Vert _{\infty }+\left\Vert
\nabla \cdot f_{2}^{2}\right\Vert _{\infty }+L_{1}^{2}+a^{+}L_{2}^{2}\right)
\int_{0}^{t}\left\Vert (Y^{1}-Y^{2})(\tau )\right\Vert _{\mathcal{H}%
}^{2}d\tau ,
\end{eqnarray*}%
which implies (\ref{700}) as claimed.

If the data are the same, this implies the uniqueness too. This ends the
proof . \hfill\ $\ $\hfill $\square $

\section{Main results}

\setcounter{equation}{0}

In this section we shall prove that the random system (\ref{26})-(\ref{29})
has a unique solution and then we shall go back through the transformation (%
\ref{23}) to the stochastic system (\ref{1})-(\ref{4}) proving that it has a
unique solution in the sense of Definition 2.1.

\medskip

\noindent \textbf{Theorem 4.1. }\textit{Under the assumptions} (\ref{31}) 
\textit{system}\textbf{\ }(\ref{26})-(\ref{29}) \textit{has, for each fixed }%
$\omega \in \Omega ,$ \textit{a unique solution }$y$\textit{, and the
process }$t\rightarrow y(t,\omega )$ \textit{is }$\mathcal{F}_{t}$\textit{%
-adapted}$.$ \textit{The solution satisfies the estimate} 
\begin{eqnarray}
&&\left\Vert y(t)\right\Vert _{\mathcal{H}}^{2}+\left\Vert y(a)\right\Vert
_{L^{2}(0,T;H)}^{2}+\int_{0}^{t}\left\Vert y(t)\right\Vert _{\mathcal{V}%
}^{2}d\tau  \label{34-1} \\
&\leq &C_{est}\left( \left\Vert y_{0}\right\Vert _{\mathcal{H}%
}^{2}+\int_{0}^{t}\left\Vert k(\tau )\right\Vert
_{L^{2}(0,a^{+};L^{2}(\partial O))}^{2}d\tau \right) ,\mbox{ \textit{for all}
}t\in \lbrack 0,T].  \notag
\end{eqnarray}

\medskip

\noindent \textbf{Proof. }We shall study first an approximating problem
introduced to endow the coefficients with more time regularity and deduce
then the necessary estimates in order to pass to the limit. Thus, we
consider a mollifier $\rho _{\varepsilon }$ and define 
\begin{equation}
W_{\varepsilon }(t,a,x)=\int_{0}^{T}W(t,a,x)\rho _{\varepsilon }(t-s)ds,
\label{34}
\end{equation}%
\begin{equation*}
\alpha _{\varepsilon }(t,a,x)=\int_{0}^{T}\alpha _{0}(t,a,x)\rho
_{\varepsilon }(t-s)ds.
\end{equation*}%
Recall that a mollifier is defined by $\rho _{\varepsilon }(t)=\frac{1}{%
\varepsilon ^{d}}\rho \left( \frac{t}{\varepsilon }\right) $ where $\rho \in
C^{\infty }(\mathbb{R}^{d}),$ $\rho (t)\geq 0,$ $\rho (t)=\rho (-t),$ $\int_{%
\mathbb{R}^{d}}\rho (t)dt=1.$ Then, $W_{\varepsilon }\in C^{\infty
}([0,T];C^{2}([0,a^{+}]\times \overline{O})),$ $\alpha _{\varepsilon }\in
C^{\infty }([0,T];L^{\infty }((0,a^{+})\times \partial O)),$ and as $%
\varepsilon \rightarrow 0$ we have\ 
\begin{eqnarray*}
W_{\varepsilon } &\rightarrow &W\mbox{ strongly in }C([0,T];C^{2}([0,a^{+}]%
\times \overline{O})), \\
\alpha _{\varepsilon } &\rightarrow &\alpha \mbox{ strongly in }L^{\infty
}((0,T)\times (0,a^{+})\times \partial O)).
\end{eqnarray*}%
The approximating system reads 
\begin{equation}
y_{t}+y_{a}-\Delta y+g_{1\varepsilon }(t,a,x)y+g_{2\varepsilon }(t,a,x)\cdot
\nabla y+\mu _{S}(t,a,x;U(e^{W_{\varepsilon }}y))y=0,\mbox{ in }(0,T)\times
(0,a^{+})\times O,  \label{35}
\end{equation}%
\begin{equation}
-\nabla y\cdot \nu =\alpha _{\varepsilon }(t,a,x)y+k(t,a,x),\mbox{ in }%
(0,T)\times (0,a^{+})\times \partial O,  \label{36}
\end{equation}%
\begin{equation}
y(t,0,x)=\int_{0}^{a^{+}}m_{\varepsilon }(t,a,x;U(e^{W_{\varepsilon
}}y))y(t,a,x)da,\mbox{ in }(0,T)\times O,  \label{37}
\end{equation}%
\begin{equation}
y(0,a,x)=y_{0}(a,x),\mbox{ in }(0,a^{+})\times O,  \label{38}
\end{equation}%
where $g_{1\varepsilon },$ $g_{2\varepsilon },$ are given by (\ref{30}) in
which $W$ is replaced by $W_{\varepsilon },$ and $\alpha $ by $\alpha
_{\varepsilon }.$ Relations (\ref{30}) and (\ref{31}) imply%
\begin{eqnarray}
g_{1\varepsilon } &\in &C^{\infty }([0,T];C^{1}([0,a^{+}]\times C(\overline{O%
})),\mbox{ \ }g_{2\varepsilon }\in C^{\infty }([0,T];C^{2}[0,a^{+}]\times
C^{1}(\overline{O})),  \label{39} \\
\alpha _{\varepsilon } &\in &C^{\infty }([0,T];L^{\infty }(0,a^{+};L^{\infty
}(\partial O))),\mbox{ \ }k\in L^{2}(0,T;L^{2}((0,a^{+})\times O).\mbox{ \ }
\notag
\end{eqnarray}%
Recall that $\mu _{S}$ and $m$ are local Lipschitz continuous with
constants, $L_{\mu _{S}}(R),$ $L_{m}(R).$

A solution to (\ref{35})-(\ref{38}) is defined by replacing in (\ref{33})
the corresponding previous coefficients.

Let us denote 
\begin{equation*}
S_{1}(t,a,x;u)=\mu _{S}(t,a,x;U(e^{W_{\varepsilon }}u))u,\mbox{ \ }%
S_{2}(t,a,x;u)=m_{0}(t,a,x;U(e^{W_{\varepsilon }}u))u,\mbox{ for }u\in 
\mathcal{H}.
\end{equation*}%
Under the local Lipschitz conditions on $m$ and $\mu _{S}$ it follows that $%
S_{1}$ and $S_{2}$ turn out to be only local Lipschitz on $\mathcal{H}.$
Indeed, let us take $R>0$ and $u,$ $\overline{u}\in \mathcal{H}$, such that $%
\left\Vert u\right\Vert _{\mathcal{H}}\leq R$ and $\left\Vert \overline{u}%
\right\Vert _{\mathcal{H}}\leq R$ and calculate 
\begin{eqnarray}
\left\vert U(e^{W_{\varepsilon }}u)\right\vert &=&\left\vert
\int_{0}^{a^{+}}\int_{O_{U}}\gamma (a,x)e^{W_{\varepsilon
}}u(a,x)dxda\right\vert \leq e^{\left\Vert W_{\varepsilon }\right\Vert
_{\infty }}\gamma _{\infty }\sqrt{a^{+}meas(O_{U})}\left\Vert u\right\Vert _{%
\mathcal{H}}  \label{81-0} \\
&\leq &c_{W}\gamma _{\infty }\sqrt{a^{+}meas(O_{U})}R,  \notag
\end{eqnarray}%
where 
\begin{equation*}
c_{W}=e^{\sup\nolimits_{t\in \lbrack 0,T]}\left\Vert W_{\varepsilon
}(t)\right\Vert _{\infty }}
\end{equation*}%
and 
\begin{eqnarray*}
&&\left\vert S_{2}(t,a,x;u)-S_{2}(t,a,x;\overline{u})\right\vert =\left\vert
m_{\varepsilon }(t,a,x;U(e^{W_{\varepsilon }}u))u-m_{\varepsilon
}(t,a,x;U(e^{W_{\varepsilon }}\overline{u}))\overline{u}\right\vert \\
&=&c_{W_{0}}\left( \left\vert (m_{0}(a,x;U(e^{W_{\varepsilon
}}u))-m_{0}(a,x;U(e^{W_{\varepsilon }}\overline{u})))\right\vert \left\vert
u\right\vert +\left\vert (u-\overline{u})m_{0}(a,x;U(e^{W_{\varepsilon }}%
\overline{u}))\right\vert \right) \\
&\leq &c_{W_{0}}c_{W}L_{m_{0}}(R)\gamma _{\infty }\sqrt{a^{+}meas(O_{U})}%
\left\Vert u-\overline{u}\right\Vert _{\mathcal{H}}\left\vert u\right\vert
+c_{W_{0}}m_{0\infty }\left\vert u-\overline{u}\right\vert
\end{eqnarray*}%
whence, denoting $C_{m}(R)=c_{W_{0}}c_{W}L_{m_{0}}(R)\gamma _{\infty }\sqrt{%
a^{+}meas(O_{U})}R+c_{W_{0}}m_{0\infty },$ we have 
\begin{equation*}
\left\Vert S_{2}(t,\cdot ,\cdot ;u)-S_{2}(\cdot ,\cdot ,\cdot ;\overline{u}%
)\right\Vert _{\mathcal{H}}\leq C_{m}(R)\left\Vert u-\overline{u}\right\Vert
_{\mathcal{H}}.
\end{equation*}%
This shows that $S_{2}$ is locally Lipschitz on $\mathcal{H}.$ A similar
relation follows for $S_{1}$ with the constant denoted $C_{\mu
_{S}}(R)=c_{W}L_{\mu _{S}}(R)\gamma _{\infty }\sqrt{a^{+}meas(O_{U})}R+\mu
_{\infty }.$

The proof will be done in two steps, the first for proving the existence of
the approximating solution and the second for passing to the limit.

\noindent \textbf{Step 1. }Let $N\geq 1.$ We approximate $S_{1}$ and $S_{2}$
by 
\begin{equation*}
S_{i}^{N}(t,a,x;u)=\left\{ 
\begin{array}{l}
S_{i}(t,a,x;u),\mbox{ \ \ \ \ \ \ \ \ }\left\Vert u\right\Vert _{\mathcal{H}%
}\leq N \\ 
S_{i}\left( t,a,x;\frac{Nu}{\left\Vert u\right\Vert _{H}}\right) ,\mbox{ \ }%
\left\Vert u\right\Vert _{\mathcal{H}}>N%
\end{array}%
\right.
\end{equation*}%
for $i=1,2.$ Then, it can be easily checked that $S_{i}^{N}(t,a,x;u)$ are
Lipschitz continuous on $\mathcal{H}$ with the constants $3C_{\mu ^{S}}(N)$
and $3C_{m}(N),$ corresponding to $i=1,2,$ respectively.

Now, we consider system (\ref{35})-(\ref{38}) with $S_{1}^{N}(t,a,x;y)$ and $%
S_{2}^{N}(t,a,x,;y)$ instead $S_{1}(t,a,x;y)$ and $S_{2}(t,a,x;y).$ In fact
this is (\ref{44})-(\ref{47}) with $E_{i}(t,a,x;y)=S_{i}^{N}(t,a,x;y)$, $%
i=1,2,$ and 
\begin{equation*}
f_{1}=g_{1\varepsilon },\mbox{ }f_{2}=g_{2\varepsilon },\mbox{ }f_{\Gamma
}=\alpha _{\varepsilon },\mbox{ }f_{\Gamma }^{0}=k,\mbox{ }Y_{0}=y_{0},\mbox{
}f=0.
\end{equation*}%
According to Proposition 3.2, this system has a unique solution $%
y_{\varepsilon }^{N}\in C([0,T];\mathcal{H})\cap
C([0,a^{+}];L^{2}(0,T;H))\cap L^{2}(0,T;\mathcal{V})$ verifying (\ref{48-0}%
), 
\begin{eqnarray}
&&-\int_{0}^{T}\int_{0}^{a^{+}}\int_{O}y_{\varepsilon }^{N}\psi
_{t}dxdadt-\int_{0}^{a^{+}}\int_{O}y_{0}\psi (0,a,x)dxda  \label{48-10} \\
&&+\int_{0}^{T}\int_{O}y_{\varepsilon }^{N}(t,a^{+},x)\psi
(t,a^{+},x)dxdt-\int_{0}^{T}\int_{O}\left(
\int_{0}^{a^{+}}S_{2}^{N}(t,a,x;y_{\varepsilon })da\right) \psi (t,0,x)dxdt 
\notag \\
&&-\int_{0}^{T}\int_{0}^{a^{+}}\int_{O}y_{\varepsilon }^{N}\psi
_{a}dxdadt+\int_{0}^{T}\int_{0}^{a^{+}}\int_{O}(\nabla y_{\varepsilon
}^{N}\cdot \nabla \psi +g_{1\varepsilon }y_{\varepsilon }^{N}\psi
+g_{2\varepsilon }\psi \cdot \nabla y_{\varepsilon }^{N})dxdadt+  \notag \\
&&\int_{0}^{T}\int_{0}^{a^{+}}\int_{\partial O}(\alpha _{\varepsilon
}y_{\varepsilon }^{N}+k)\psi d\sigma
dadt+\int_{0}^{T}\int_{0}^{a^{+}}\int_{O}S_{1}^{N}(t,a,x;y_{\varepsilon
})\psi dxdadt=0.\mbox{ }  \notag
\end{eqnarray}%
Moreover, the solution satisfies the estimates (\ref{48-1}), 
\begin{eqnarray}
&&\left\Vert y_{\varepsilon }^{N}(t)\right\Vert _{\mathcal{H}%
}^{2}+\left\Vert y_{\varepsilon }^{N}(a)\right\Vert
_{L^{2}(0,T;H)}^{2}+\int_{0}^{t}\left\Vert y_{\varepsilon
}^{N}(t)\right\Vert _{\mathcal{V}}^{2}d\tau  \label{507} \\
&\leq &C_{est}\left( \left\Vert y_{0}\right\Vert _{\mathcal{H}%
}^{2}+\int_{0}^{t}\left\Vert k(\tau )\right\Vert
_{L^{2}(0,a^{+};L^{2}(\partial O))}^{2}d\tau \right)  \notag
\end{eqnarray}%
and two solutions corresponding to two sets of data obey the inequality (\ref%
{700}).

Here, $C_{est}=c_{0}e^{c_{1}(1+\left\Vert g_{1\varepsilon }\right\Vert
_{\infty }+\left\Vert g_{2\varepsilon }\right\Vert _{\infty
}^{2}+a^{+}m_{\infty }^{2}+\mu _{\infty })},$ by (\ref{Cest}), where $%
\left\Vert g_{i\varepsilon }\right\Vert _{\infty }\leq \left\Vert
g_{i}\right\Vert _{\infty }\leq C_{i}$ (depending on $\left\Vert
W_{a}\right\Vert _{\infty },$ $\left\Vert \Delta W\right\Vert _{\infty })$ , 
$i=1,2,$ because the functions $g_{i}$ are continuous$.$

Now, we set 
\begin{equation}
R_{0}:=c_{0}e^{c_{1}(1+\left\Vert g_{1}\right\Vert _{\infty }+\left\Vert
g_{2}\right\Vert _{\infty }^{2}+a^{+}m_{\infty }^{2}+\mu _{\infty })}\left(
\left\Vert y_{0}\right\Vert _{\mathcal{H}}^{2}+\int_{0}^{T}\left\Vert
k(t)\right\Vert _{L^{2}(0,a^{+};L^{2}(\partial O))}^{2}dt\right) .
\label{R0}
\end{equation}%
It follows that for $N\geq \lbrack R_{0}]+1:=N_{0}$ we get $\left\Vert
y_{\varepsilon }^{N}(t)\right\Vert _{\mathcal{H}}^{2}\leq R_{0}<N_{0}\leq N$
and so, $S_{i}^{N}(t,a,x;y)=S_{i}(t,a,x;y)$, $i=1.2,$ meaning that $%
y_{\varepsilon }^{N}$ actually satisfies system (\ref{35})-(\ref{38}), if $%
N\geq N_{0}.$ Thus, we deduce that $y_{\varepsilon }^{N_{0}}(t)$ is in fact
a solution to problem (\ref{35})-(\ref{38}) and we denote it by $%
y_{\varepsilon }(t).$ We also note that the Lipschitz constants for $%
S_{i}^{N}$ specified before depend actually on $R_{0},$ namely%
\begin{eqnarray*}
L_{1} &=&C_{m}(R_{0})=c_{W_{0}}c_{W}L_{m_{0}}(R_{0})\gamma _{\infty }\sqrt{%
a^{+}meas(O_{U})}R+c_{W_{0}}m_{0\infty },\mbox{ \ } \\
L_{2} &=&C_{\mu _{S}}(R_{0})=c_{W}L_{\mu _{S}}(R_{0})\gamma _{\infty }\sqrt{%
a^{+}meas(O_{U})}R+\mu _{\infty }.
\end{eqnarray*}

To prove the uniqueness, we consider two solutions $y_{\varepsilon }$ and $%
\overline{y_{\varepsilon }}$ corresponding to the same data and see that for 
$N>\sup\limits_{t\in \lbrack 0,T]}\left( \left\Vert y_{\varepsilon
}^{N}(t)\right\Vert _{\mathcal{H}}^{2}+\left\Vert \overline{y_{\varepsilon
}^{N}}(t)\right\Vert _{\mathcal{H}}^{2}\right) $ it follows by (\ref{700})
that their difference is zero.

Obviously, the solution $y_{\varepsilon }^{N_{0}}=y_{\varepsilon }$
satisfies (\ref{48-t}), in which $E_{i}$ are replaced by $%
S_{i}(t,a,x;U(e^{W_{\varepsilon }}y_{\varepsilon }))$, $i=1,2,$ 
\begin{eqnarray}
&&\int_{0}^{a}\int_{O}y_{\varepsilon }(t,s,x)\psi
(t,s,x)dxda-\int_{0}^{t}\int_{0}^{a}\int_{O}y_{\varepsilon }\psi _{\tau
}dxdsd\tau -\int_{0}^{a}\int_{O}y_{0}\psi (0,s,x)dxds  \label{509} \\
&&+\int_{0}^{t}\int_{O}y_{\varepsilon }(\tau ,a,x)\psi (\tau ,a,x)dxd\tau
-\int_{0}^{t}\int_{0}^{a}\int_{O}y_{\varepsilon }\psi _{a}dxdsd\tau  \notag
\\
&&-\int_{0}^{t}\int_{O}\left( \int_{0}^{a}m_{\varepsilon }(\tau
,s,x;U(e^{W_{\varepsilon }}y_{\varepsilon }))y_{\varepsilon }ds\right) \psi
(\tau ,0,x)dxd\tau  \notag \\
&&+\int_{0}^{t}\int_{0}^{a}\int_{O}(\nabla y_{\varepsilon }\cdot \nabla \psi
+g_{1\varepsilon }y_{\varepsilon }\psi +g_{2\varepsilon }\psi \cdot \nabla
y_{\varepsilon })dxdsd\tau +  \notag \\
&&\int_{0}^{t}\int_{0}^{a}\int_{\partial O}(\alpha _{\varepsilon
}y_{\varepsilon }\psi +k_{\varepsilon }\psi )d\sigma dsd\tau
+\int_{0}^{t}\int_{0}^{a}\int_{O}\mu _{S}(\tau ,s,x;U(e^{W_{\varepsilon
}}y_{\varepsilon }))y_{\varepsilon }\psi dxdsd\tau =0,  \notag
\end{eqnarray}%
and inherits estimates (\ref{507}). Moreover, (\ref{600}), written for $%
Y_{1}=y_{\varepsilon },$ $Y_{2}=y_{\varepsilon ^{\prime }},$ corresponding
to $W_{\varepsilon }$ and $W_{\varepsilon ^{\prime }},$ respectively, yields%
\begin{eqnarray}
&&\left\Vert (y_{\varepsilon }-y_{\varepsilon ^{\prime }})(t)\right\Vert _{%
\mathcal{H}}^{2}+\int_{0}^{t}\int_{O}(y_{\varepsilon }-y_{\varepsilon
^{\prime }})^{2}(\tau ,a,x)dxd\tau +\int_{0}^{t}\left\Vert (y_{\varepsilon
}-y_{\varepsilon ^{\prime }})(\tau )\right\Vert _{\mathcal{V}}^{2}d\tau
\label{509-1} \\
&\leq &c_{0}R_{0}^{2}\left( \left\Vert g_{1\varepsilon }-g_{1\varepsilon
^{\prime }}\right\Vert _{\infty }^{2}+\left\Vert g_{2\varepsilon
}-g_{2\varepsilon ^{\prime }}\right\Vert _{\infty }^{2}+\left\Vert \alpha
_{\varepsilon }-\alpha _{\varepsilon ^{\prime }}\right\Vert _{\infty
}^{2}\right)  \notag \\
&&+c_{1}\left( 1+\left\Vert g_{1\varepsilon ^{\prime }}\right\Vert _{\infty
}+\left\Vert \nabla \cdot g_{2\varepsilon ^{\prime }}\right\Vert _{\infty
}\right) \int_{0}^{t}\left\Vert (y_{\varepsilon }-y_{\varepsilon ^{\prime
}})(\tau )\right\Vert _{\mathcal{H}}^{2}d\tau  \notag \\
&&+c_{1}\int_{0}^{t}\int_{O}\left( \int_{0}^{a^{+}}(S_{2}(\tau
,a,x;y_{\varepsilon })-S_{2}(\tau ,a,x;y_{\varepsilon ^{\prime }}))ds\right)
^{2}dxd\tau  \notag \\
&&+c_{1}\int_{0}^{t}\int_{0}^{a^{+}}\int_{O}(S_{1}(\tau ,a,x;y_{\varepsilon
})-S_{1}(\tau ,a,x;y_{\varepsilon ^{\prime }}))^{2}dxdsd\tau .  \notag
\end{eqnarray}%
\noindent \textbf{Step 2.} The second step is devoting to passing to the
limit as $\varepsilon \rightarrow 0.$ To this end, we use (\ref{509-1}) and
detail first some computations.

Recall that by (\ref{81-0}), $\left\vert U(e^{W_{\varepsilon
}}y_{\varepsilon })\right\vert \leq c_{W}\gamma _{\infty }\sqrt{%
a^{+}meas(O_{U})}R_{0},$ where $R_{0}$ is precisely (\ref{R0}), and we
calculate 
\begin{eqnarray*}
&&\left\vert S_{2}(t,a,x;y_{\varepsilon })-S_{2}(t,a,x;y_{\varepsilon
^{\prime }})\right\vert =\left\vert m_{\varepsilon
}(t,a,x;U(e^{W_{\varepsilon }}y_{\varepsilon }))y_{\varepsilon
}-m_{\varepsilon ^{\prime }}(t,a,x;U(e^{W_{\varepsilon }}y_{\varepsilon
^{\prime }}))y_{\varepsilon ^{\prime }}\right\vert \\
&=&\left\vert (m_{\varepsilon }(t,a,x;U(e^{W_{\varepsilon }}y_{\varepsilon
}))-m_{\varepsilon ^{\prime }}(t,a,x;U(e^{W_{\varepsilon }}y_{\varepsilon
^{\prime }}))\right\vert \left\vert y_{\varepsilon }\right\vert +\left\vert
(y_{\varepsilon }-y_{\varepsilon ^{\prime }})m_{\varepsilon ^{\prime
}}(t,a,x;U(e^{W_{\varepsilon }}y_{\varepsilon ^{\prime }}))\right\vert .
\end{eqnarray*}%
Recall that $\left\{ e^{W_{\varepsilon }(t,a,x)-W_{\varepsilon
}(t,0,x)}\right\} _{\varepsilon }$ is a Cauchy sequence and by (\ref{30}) we
have 
\begin{eqnarray*}
&&\left\vert (m_{\varepsilon }(t,a,x;U(e^{W_{\varepsilon }}y_{\varepsilon
}))-m_{\varepsilon ^{\prime }}(t,a,x;U(e^{W_{\varepsilon }}y_{\varepsilon
^{\prime }}))\right\vert \\
&=&\left\vert (m_{0}(a,x;U(e^{W_{\varepsilon }}y_{\varepsilon
}))e^{W_{\varepsilon }(t,a,x)-W_{\varepsilon
}(t,0,x)}-m_{0}(a,x;U(e^{W_{\varepsilon ^{\prime }}}y_{\varepsilon ^{\prime
}}))e^{W_{\varepsilon ^{\prime }}(t,a,x)-W_{\varepsilon ^{\prime
}}(t,0,x)}\right\vert \\
&=&\left\vert (m_{0}(a,x;U(e^{W_{\varepsilon }}y_{\varepsilon
}))-m_{0}(a,x;U(e^{W_{\varepsilon ^{\prime }}}y_{\varepsilon ^{\prime
}}))\right\vert \left\vert e^{W_{\varepsilon }(t,a,x)-W_{\varepsilon
}(t,0,x)}\right\vert \\
&&+\left\vert e^{W_{\varepsilon }(t,a,x)-W_{\varepsilon
}(t,0,x)}-e^{W_{\varepsilon ^{\prime }}(t,a,x)-W_{\varepsilon ^{\prime
}}(t,0,x)}\right\vert \left\vert m_{0}(a,x;U(e^{W_{\varepsilon ^{\prime
}}}y_{\varepsilon ^{\prime }}))\right\vert \\
&\leq &c_{W_{0}}L_{m_{0}}(R_{0})\gamma _{\infty
}\int_{0}^{a^{+}}\int_{O_{U}}\left( \left\vert e^{W_{\varepsilon
}}-e^{W_{\varepsilon ^{\prime }}}\right\vert \left\vert y_{\varepsilon
}(t)\right\vert +\left\vert y_{\varepsilon }-y_{\varepsilon ^{\prime
}}\right\vert \left\vert e^{W_{\varepsilon ^{\prime }}}\right\vert \right)
dxda+m_{0\infty }\delta _{\varepsilon ,\varepsilon ^{\prime }}
\end{eqnarray*}%
with $\delta _{\varepsilon ,\varepsilon ^{\prime }}$ arbitrarily small. Then,%
\begin{eqnarray*}
&&\left\vert S_{2}(t,a,x;y_{\varepsilon }(t))-S_{2}(t,a,x;y_{\varepsilon
^{\prime }}(t))\right\vert \\
&\leq &(c_{W_{0}}L_{m_{0}}(R)\gamma _{\infty }\sqrt{a^{+}meas(O_{U})}%
\left\Vert y_{\varepsilon }(t)\right\Vert _{\mathcal{H}}\delta _{\varepsilon
,\varepsilon ^{\prime }}+\left\Vert y_{\varepsilon }-y_{\varepsilon ^{\prime
}}\right\Vert _{\mathcal{H}}c_{W}+m_{0\infty }\delta _{\varepsilon
,\varepsilon ^{\prime }})\left\vert y_{\varepsilon }(t)\right\vert _{%
\mathcal{H}}+m_{\infty }\left\vert y_{\varepsilon }-y_{\varepsilon ^{\prime
}}\right\vert
\end{eqnarray*}%
whence 
\begin{equation*}
\left\Vert S_{2}(\cdot ,\cdot ,\cdot ;y_{\varepsilon }(t))-S_{2}(\cdot
,\cdot ,\cdot ;y_{\varepsilon ^{\prime }}(t))\right\Vert _{\mathcal{H}}\leq
C_{m}(R_{0})\left\Vert y_{\varepsilon }(t)-y_{\varepsilon ^{\prime
}}(t)\right\Vert _{\mathcal{H}}+C_{2}(R_{0})\delta _{\varepsilon
,\varepsilon ^{\prime }}.
\end{equation*}%
For $S_{1}$ we get%
\begin{equation*}
\left\Vert S_{1}(\cdot ,\cdot ,\cdot ;y_{\varepsilon }(t))-S_{1}(\cdot
,\cdot ,\cdot ;y_{\varepsilon ^{\prime }}(t))\right\Vert _{\mathcal{H}}\leq
C_{\mu _{S}}(R_{0})\left\Vert y_{\varepsilon }(t)-y_{\varepsilon ^{\prime
}}(t)\right\Vert _{\mathcal{H}}+C_{3}(R_{0})\delta _{\varepsilon
,\varepsilon ^{\prime }}.
\end{equation*}%
Then, (\ref{509-1}) yields 
\begin{eqnarray*}
&&\left\Vert (y_{\varepsilon }-y_{\varepsilon ^{\prime }})(t)\right\Vert _{%
\mathcal{H}}^{2}+\int_{0}^{t}\int_{O}(y_{\varepsilon }-y_{\varepsilon
^{\prime }})^{2}(\tau ,a,x)dxd\tau +\int_{0}^{t}\left\Vert (y_{\varepsilon
}-y_{\varepsilon ^{\prime }})(\tau )\right\Vert _{\mathcal{V}}^{2}d\tau \\
&\leq &c_{0}R_{0}^{2}\left( \left\Vert g_{1\varepsilon }-g_{1\varepsilon
^{\prime }}\right\Vert _{\infty }^{2}+\left\Vert g_{2\varepsilon
}-g_{2\varepsilon ^{\prime }}\right\Vert _{\infty }^{2}+\left\Vert \alpha
_{\varepsilon }-\alpha _{\varepsilon ^{\prime }}\right\Vert _{\infty
}^{2}\right) \\
&&+c_{1}\left( 1+\left\Vert g_{1\varepsilon ^{\prime }}\right\Vert _{\infty
}+\left\Vert \nabla \cdot g_{2\varepsilon ^{\prime }}\right\Vert _{\infty
}\right) \int_{0}^{t}\left\Vert (y_{\varepsilon }-y_{\varepsilon ^{\prime
}})(\tau )\right\Vert _{\mathcal{H}}^{2}d\tau \\
&&+c_{1}\left( (a^{+}C_{m}^{2}(R_{0})+C_{\mu
_{S}}^{2}(R_{0}))\int_{0}^{t}\left\Vert y_{\varepsilon }(\tau
)-y_{\varepsilon ^{\prime }}(\tau )\right\Vert _{\mathcal{H}}^{2}d\tau
+(C_{2}^{2}(R_{0})+C_{3}^{2}(R_{0}))\delta _{\varepsilon ,\varepsilon
^{\prime }}^{2}\right) ,
\end{eqnarray*}%
and applying the Gronwall's lemma we get 
\begin{eqnarray*}
&&\left\Vert (y_{\varepsilon }-y_{\varepsilon ^{\prime }})(t)\right\Vert _{%
\mathcal{H}}^{2}+\int_{0}^{t}\int_{O}(y_{\varepsilon }-y_{\varepsilon
^{\prime }})^{2}(\tau ,a,x)dxd\tau +\int_{0}^{t}\left\Vert (y_{\varepsilon
}-y_{\varepsilon ^{\prime }})(\tau )\right\Vert _{\mathcal{V}}^{2}d\tau \\
&\leq &e^{c_{1}\left( 1+\left\Vert g_{1\varepsilon ^{\prime }}\right\Vert
_{\infty }+\left\Vert \nabla \cdot g_{2\varepsilon ^{\prime }}\right\Vert
_{\infty }+(a^{+}C_{m}^{2}(R_{0})+C_{\mu _{S}}^{2}(R_{0}))\right) T} \\
&&\times c_{0}R_{0}^{2}\left( \left\Vert g_{1\varepsilon }-g_{1\varepsilon
^{\prime }}\right\Vert _{\infty }^{2}+\left\Vert g_{2\varepsilon
}-g_{2\varepsilon ^{\prime }}\right\Vert _{\infty }^{2}+\left\Vert \alpha
_{\varepsilon }-\alpha _{\varepsilon ^{\prime }}\right\Vert _{\infty
}^{2}+c_{1}(C_{2}^{2}(R_{0})+C_{3}^{2}(R_{0}))\delta _{\varepsilon
,\varepsilon ^{\prime }}^{2}\right) .
\end{eqnarray*}%
Taking into account that $\{W_{\varepsilon }\}_{\varepsilon },$ $\{\alpha
_{\varepsilon }\}_{\varepsilon },$ are Cauchy sequences, we deduce that $%
\{y_{\varepsilon }\}_{\varepsilon >0}$ is a Cauchy sequence too, hence 
\begin{equation*}
y_{\varepsilon }\rightarrow y\mbox{ strongly in }C([0,T];\mathcal{H})\cap
L^{2}(0,T;\mathcal{V})\cap C([0,a^{+}];L^{2}(0,T;H)).
\end{equation*}%
Consequently, since%
\begin{equation*}
\left\vert U(y_{\varepsilon })-U(y)\right\vert =\left\vert
\int_{0}^{a^{+}}\int_{O_{U}}\gamma (a,x)(y_{\varepsilon
}-y)(t)dxda\right\vert \leq \gamma _{\infty }\sqrt{a^{+}meas(O_{U})}%
\left\Vert y_{\varepsilon }(t)-y(t)\right\Vert _{\mathcal{H}},
\end{equation*}%
it follows that $U(y_{\varepsilon })\rightarrow U(y)$ strongly in $C([0,T]),$
and a.e. on $(0,T).$ By Egorov theorem, there exists a measurable subset $%
A_{\delta ^{\prime }}\subset (0,T),$ with $meas(A_{\delta ^{\prime
}})<\delta ^{\prime },$ and $U(y_{\varepsilon })\rightarrow U(y)$ uniformly
on $(0,T)\backslash A_{\delta ^{\prime }}.$ Then, since $m_{0}$ is
continuous with respect to the fourth variable we have 
\begin{equation*}
m_{\varepsilon }(t,a,x;U(y_{\varepsilon }))=e^{W_{\varepsilon
}(t,a,x)-W_{\varepsilon }(t,0,x)}m_{0}(a,x;U(y_{\varepsilon }))\rightarrow
m_{0}(a,x;U(y))e^{W(t,a,x)-W(t,0,x)}
\end{equation*}%
on $(0,T)\backslash A_{\delta ^{\prime }}\times (0,a^{+})\times O$ and so it
tends strongly in $L^{2}(0,T;\mathcal{H}).$ This implies 
\begin{equation*}
m_{\varepsilon }(\cdot ,\cdot ,\cdot ;U(y_{\varepsilon }))y_{\varepsilon
}\rightarrow m(\cdot ,\cdot ,\cdot ;U(y))y\mbox{ strongly in }L^{2}(0,T;%
\mathcal{H}).
\end{equation*}%
A similar convergence is true for $\mu _{S}(\cdot ,\cdot ,\cdot ;U(t,y))y.$
Since the coefficients $g_{1\varepsilon },$ $g_{2\varepsilon },$ $\alpha
_{\varepsilon },$ $k_{\varepsilon }$ tend strongly to $g_{1},$ $g_{2},$ $%
\alpha ,$ $k$ in their corresponding spaces, it follows by passing to the
limit in (\ref{509}) that $y$ satisfies (\ref{33}). In particular for $t=T,$ 
$a=a^{+}$ $\psi (T,a,x)=0,$ it is in conclusion a solution to (\ref{26})-(%
\ref{29}).

Relations (\ref{507}) and (\ref{700}) are satisfied at limit by $y$ and the
difference $y-\overline{y},$ respectively, and imply (\ref{34-1}) and the
solution uniqueness.

Next, we show that $y(t)$ is a $\mathcal{F}_{t}$-adapted process. Let us
recall problem (\ref{49}) and assume that $Y_{0}$ is measurable with respect
to $\mathcal{F}_{0}$. In fact, $Y_{0}$ stands for $y_{0}=p_{0}$ which has
this property by (\ref{11-0}). Since $A(t)$ is quasi $m$-accretive, one can
consider this equation with $A(t)$ replaced by its Yosida approximation $%
A_{\lambda }(t)$ which is Lipschitz. The solution to the approximating
equation can be obtained by an iterative process and so it is measurable
with respect to $\mathcal{F}_{t}$. Also this property is preserved by
passing to the limit, then the solution to (\ref{49}), as well as all the
other solutions, that is $v_{\varepsilon }(t),$ $y_{\varepsilon }(t)$ and $%
y(t)$ in Theorem 3.1 which are deduced as limits of $\mathcal{F}_{t}$%
-adapted sequences, so that they are $\mathcal{F}_{t}$-adapted$.$

The proof is ended.\hfill $\square $

\medskip

In addition to the properties of $y$ proved in Theorem 3.1 one can add that,
for each $\omega \in \Omega ,$ there exists the strong derivative of $y$ and
equations (\ref{26})-(\ref{29}) are satisfied in the sense of distributions$%
. $

Let us define 
\begin{equation*}
\mathcal{X}=\{u\in \mathcal{V};\mbox{ }u_{a}\in \mathcal{V}^{\prime }\},%
\mbox{ }H_{T}^{1}(0,T)=\left\{ \varphi \in H^{1}(0,T);\varphi (T)=0\right\} 
\mbox{ }
\end{equation*}%
and denote by $\mathcal{X}^{\prime }$ and $(H_{T}^{1}(0,T))^{\prime }$ the
dual spaces of $\mathcal{X}$ and $H_{T}^{1}(0,T),$ respectively.

\medskip

\noindent \textbf{Corollary 4.2. }\textit{Under the assumptions} \textit{of}
Theorem 4.1 \textit{it follows that} 
\begin{equation}
\frac{dy}{dt}\in L^{2}(0,T;\mathcal{X}^{\prime }).  \label{94-00}
\end{equation}

\medskip

\noindent \textbf{Proof}. In (\ref{48-0}) $\psi $ can be taken of the form $%
\psi (t,a,x)=\varphi (t)\psi _{0}(a,x),$ with $\varphi \in H_{T}^{1}(0,T)$
and $\psi _{0}\in \mathcal{X}.$ Obviously, $\psi _{0}\in C([0,a^{+}];H).$
Let us define $\widetilde{A}(t):\mathcal{V}\cap C[0,a^{+};H]\rightarrow 
\mathcal{X}^{\prime }$ by 
\begin{eqnarray}
&&\left\langle \widetilde{A}(t)v,\psi _{0}\right\rangle _{\mathcal{X}%
^{\prime },\mathcal{X}}=\int_{O}v(a^{+},x)\psi
_{0}(a^{+},x)dx-\int_{0}^{a^{+}}\int_{O}v(\psi _{0})_{a}dxda  \label{93-0} \\
&&-\int_{O}\left( \int_{0}^{a^{+}}m(t,a,x;U(v))vda\right) \psi
_{0}(0,x)dx+\int_{0}^{a^{+}}\int_{O}\mu _{S}(t,a,x;U(v))v\psi _{0}dxda 
\notag \\
&&+\int_{0}^{a^{+}}\int_{O}(\nabla v\cdot \nabla \psi _{0}+vg_{1}\psi
_{0}+\psi _{0}g_{2}\cdot \nabla v)dxda+\int_{0}^{a^{+}}\int_{\partial
O}(\alpha v+k)\psi _{0}d\sigma da,\mbox{ }  \notag
\end{eqnarray}%
for all $v\in \mathcal{X},$ where $\mathcal{X}^{\prime }$ is the dual of $%
\mathcal{X},$ with the pivot space $\mathcal{H}.$

One can easily calculate that $\left\Vert \widetilde{A}(t)v\right\Vert _{%
\mathcal{X}^{\prime }}\leq C\left( \left\Vert v\right\Vert _{\mathcal{V}%
}+\left\Vert v\right\Vert _{C([0,a^{+}];H)}\right) ,$ hence $\widetilde{A}%
(t) $ is well defined$.$

Moreover, for any $\varphi \in H_{T}^{1}(0,T)$ and $\psi _{0}\in \mathcal{X}$%
, we define the distributional derivative%
\begin{equation}
\frac{dy}{dt}(\varphi )=-\int_{0}^{T}y(t,a,x)\varphi
_{t}(t)dt-y_{0}(a,x)\varphi (0),  \label{93}
\end{equation}%
and%
\begin{equation*}
\left\langle \frac{dy}{dt}(\varphi ),\psi _{0}\right\rangle _{\mathcal{X}%
^{\prime },\mathcal{X}}=-\int_{0}^{T}\int_{0}^{a^{+}}\int_{O}y(t,a,x)\varphi
_{t}(t)\psi _{0}dxdadt-\int_{0}^{a^{+}}\int_{O}y_{0}(a,x)\varphi (0)\psi
_{0}dxda.
\end{equation*}%
Then, one can write (\ref{48-0}) as 
\begin{equation}
\left\langle \frac{dy}{dt}(\varphi ),\psi _{0}\right\rangle _{\mathcal{X}%
^{\prime },\mathcal{X}}+\int_{0}^{T}\left\langle \widetilde{A}%
(t)y(t),\varphi (t)\psi _{0}\right\rangle _{\mathcal{X}^{\prime },\mathcal{X}%
}dt=0,  \label{92}
\end{equation}%
for any $\varphi \in H_{T}^{1}(0,T)$ and $\psi _{0}\in \mathcal{X}.$ This
implies%
\begin{equation*}
\frac{dy}{dt}(\varphi )+\widetilde{A}(t)y(\varphi )=0,\mbox{ for all }%
\varphi \in H_{T}^{1}(0,T),
\end{equation*}%
which can be still written%
\begin{equation}
\frac{dy}{dt}+\widetilde{A}(t)y=0,\mbox{ in }\mathcal{D}^{\prime }(0,T;%
\mathcal{X}^{\prime }).  \label{94}
\end{equation}%
Moreover, since 
\begin{eqnarray*}
&&\int_{0}^{T}\left\vert \left\langle \widetilde{A}(t)y(t),\varphi (t)\psi
_{0}\right\rangle _{\mathcal{X}^{\prime },\mathcal{X}}\right\vert dt\leq
\int_{0}^{T}\left\vert \varphi (t)\right\vert \left\Vert A(t)y(t)\right\Vert
_{\mathcal{X}^{\prime }}\left\Vert \psi _{0}\right\Vert _{\mathcal{X}}dt \\
&\leq &\left\Vert \varphi \right\Vert _{L^{2}(0,T)}\left(
\int_{0}^{T}\left\Vert A(t)y(t)\right\Vert _{\mathcal{X}^{\prime
}}^{2}dt\right) ^{1/2}\left\Vert \psi _{0}\right\Vert _{\mathcal{X}}\leq
C\left\Vert \varphi \right\Vert _{L^{2}(0,T)}\left\Vert \psi _{0}\right\Vert
_{\mathcal{X}}\left( \left\Vert v\right\Vert _{\mathcal{V}}+\left\Vert
v\right\Vert _{C([0,a^{+}];H)}\right)
\end{eqnarray*}%
it follows that 
\begin{equation*}
\left\Vert \frac{dy}{dt}(\varphi )\right\Vert _{\mathcal{X}^{\prime }}\leq
\sup_{\left\Vert \psi _{0}\right\Vert _{\mathcal{X}}\leq
1}\int_{0}^{T}\left\vert \left\langle \widetilde{A}(t)y(t),\varphi (t)\psi
_{0}\right\rangle _{\mathcal{X}^{\prime },\mathcal{X}}\right\vert dt\leq
C\left\Vert \varphi \right\Vert _{L^{2}(0,T)}
\end{equation*}%
implying (\ref{94-00}), and so (\ref{94}) can be written 
\begin{equation}
\frac{dy}{dt}(t)+\widetilde{A}(t)y(t)=0,\mbox{ a.e. }t\in (0,T).
\label{94-1}
\end{equation}

\medskip

\noindent \textbf{Theorem 4.3. }\textit{Under the assumptions} (\ref{9-Lip}%
)-(\ref{6}) \textit{the stochastic problem }(\ref{1})-(\ref{4}) \textit{has
a unique solution and }$e^{-w}p\in L^{2}(0,T;\mathcal{X}^{\prime }).$

\medskip

\noindent \textbf{Proof. }Recall that (\ref{26})-(\ref{29}) has a unique
solution (\ref{33}), for each $\omega \in \Omega ,$ given by Theorem 3.1. We
go back to $p$ by the transformation (\ref{23}).

Let us consider again a mollifier $\rho _{\varepsilon }$ and define the
function%
\begin{equation*}
y_{\varepsilon }(t)=(y\ast \rho _{\varepsilon })(t)=\int_{0}^{T}y(t-s)\rho
_{\varepsilon }(s)ds.
\end{equation*}%
Obviously, 
\begin{equation}
y_{\varepsilon }\rightarrow y\mbox{ strongly in }C([0,T];\mathcal{H})\cap
C([0,a^{+}];L^{2}(0,T;H))\cap L^{2}(0,T;\mathcal{V}),  \label{95}
\end{equation}%
and $\left\Vert y_{\varepsilon }(t)\right\Vert _{\mathcal{H}}\leq \left\Vert
y(t)\right\Vert _{\mathcal{H}}$ which satisfies (\ref{34-1}).

We have 
\begin{equation*}
\frac{dy_{\varepsilon }}{dt}(t)=\frac{d}{dt}(y\ast \rho _{\varepsilon
})(t)=\int_{0}^{T}\frac{dy}{dt}(t-s)\rho _{\varepsilon }(s)ds=\left( \frac{dy%
}{dt}\ast \rho _{\varepsilon }\right) (t).
\end{equation*}%
We multiply (\ref{94-1}) by $\rho _{\varepsilon },$ and get%
\begin{equation*}
\frac{dy_{\varepsilon }}{dt}(t)+(\rho _{\varepsilon }\ast \widetilde{A}%
(t)y_{\varepsilon })(t)=0,
\end{equation*}%
and then by $e^{W},$ obtaining%
\begin{equation}
e^{W}\frac{dy_{\varepsilon }}{dt}(t)+e^{W}(\rho _{\varepsilon }\ast 
\widetilde{A}(t)y_{\varepsilon })(t)=0.  \label{96}
\end{equation}%
Let us denote $p_{\varepsilon }:=e^{W}y_{\varepsilon }$ and note that $%
p_{\varepsilon }\rightarrow e^{W}y:=p$ strongly in all spaces indicated in (%
\ref{95}).

Next, by It\^{o}'s formula we have%
\begin{equation*}
e^{W}dy_{\varepsilon }=d(e^{W}y_{\varepsilon })-y_{\varepsilon }de^{W}
\end{equation*}%
and using (\ref{23-0}) in (\ref{96}) we get%
\begin{equation*}
dp_{\varepsilon }-p_{\varepsilon }dW-\mu p_{\varepsilon }dt+e^{W}(\rho
_{\varepsilon }\ast \widetilde{A}(t)y_{\varepsilon })(t)dt=0.
\end{equation*}%
Integrating from $0$ to $t$ and taking into account that 
\begin{equation*}
\int_{0}^{T}\left( \int_{0}^{a^{+}}\int_{O}p_{\varepsilon }(t)\mu
_{j}dxda\right) ^{2}dt\leq \left\Vert \mu _{j}\right\Vert _{\infty
}\int_{0}^{T}\left\Vert p_{\varepsilon }(t)\right\Vert _{\mathcal{H}%
}^{2}dt\leq C,\mbox{ }\mathbb{P}\mbox{-a.s.,}
\end{equation*}%
which ensures that the It\^{o} integral makes sense, we have%
\begin{equation*}
p_{\varepsilon }(t)-p_{\varepsilon }(0)-\int_{0}^{t}\mu p_{\varepsilon
}d\tau -\int_{0}^{t}p_{\varepsilon }(\tau )dW(\tau )+\int_{0}^{t}e^{W(\tau
)}(\rho _{\varepsilon }\ast \widetilde{A}(\tau )y_{\varepsilon })(\tau
)d\tau =0.
\end{equation*}%
Then, passing to the limit as $\varepsilon \rightarrow 0$ and taking into
account the convergence of $p_{\varepsilon }$ to $p$ and the definition (\ref%
{93-0}) we obtain%
\begin{equation*}
p(t)-p(0)-\int_{0}^{t}\mu pd\tau -\int_{0}^{t}p(\tau )dW(\tau
)+\int_{0}^{t}e^{W(\tau ,a,x)}\widetilde{A}(\tau )y(\tau )d\tau =0.
\end{equation*}%
This equation tested at $\psi _{0}\in \mathcal{X},$ yields%
\begin{eqnarray}
&&\left( p(t),\psi _{0}\right) _{\mathcal{H}}-\left( p(0),\psi _{0}\right) _{%
\mathcal{H}}-\int_{0}^{t}\left( \mu p(\tau ),\psi _{0}\right) _{\mathcal{H}%
}d\tau -\int_{0}^{t}\left( p(\tau ),\psi _{0}\right) _{\mathcal{H}}dW(\tau )
\label{97} \\
&&+\int_{0}^{t}\left\langle e^{W(\tau )}\widetilde{A}(\tau )y(\tau ),\psi
_{0}\right\rangle _{\mathcal{X}^{\prime },\mathcal{X}}d\tau =0.  \notag
\end{eqnarray}%
But 
\begin{equation*}
\left\langle e^{W(\tau )}\widetilde{A}(\tau )y(\tau ),\psi _{0}\right\rangle
_{\mathcal{X}^{\prime },\mathcal{X}}=\left\langle \widetilde{A}(\tau )y(\tau
),e^{W(\tau )}\psi _{0}\right\rangle _{\mathcal{X}^{\prime },\mathcal{X}},
\end{equation*}%
and replacing in (\ref{97}) the definition of $\widetilde{A}(\tau )y(\tau )$
by (\ref{93-0}) where the test function is $e^{W(\tau )}\psi _{0}$ we obtain
after performing all necessary calculations the weak form (\ref{21}).

The solution $p$ is constructed as the limit of an $\mathcal{F}_{t}$-adapted
sequence, so that $p$ is a $\mathcal{F}_{t}$-adapted process.

Finally, let us assume that there are two solutions $p_{1}$ and $p_{2}$
satisfying (\ref{1})-(\ref{4}). By substituting $y_{i}=e^{-W}p_{i},$ $i=1,2,$%
and by making all calculations we are led to two systems in $y_{i}$ with the
same coefficients. As we know that the solution to the deterministic random
system (\ref{26})-(\ref{29}) is unique, it follows that the solution $p$ to
the stochastic system in unique. This ends the proof.\hfill $\square $

\bigskip

\noindent \textbf{Acknowledgement.} This work was supported by a grant of
Ministry of Research and Innovation, CNCS -- UEFISCDI, project number
PN-III-P4-ID-PCE-2016-0011, within PNCDI III.

\end{document}